\documentclass[11pt]{article}

\usepackage[utf8]{inputenc}
\usepackage[T1]{fontenc}
\usepackage{graphicx}
\usepackage{latexsym}
\usepackage{amsmath,amssymb,amsthm}
\usepackage{mathrsfs}
\usepackage{xcolor}
\usepackage{stmaryrd}
\usepackage{bbm}
\usepackage{url}
\usepackage{longtable}
\usepackage{dsfont}
\usepackage{lmodern}

\usepackage[numbers]{natbib}
\usepackage{color}
\usepackage{lmodern}
\usepackage{float}
\usepackage{enumitem}
\setenumerate{leftmargin=*}
\usepackage{amsopn}
\usepackage{comment}

\usepackage[linktocpage]{hyperref}

\newcommand{\bbr}{\mathbb{R}}

\newcommand{\eps}{{\varepsilon}}

\newcommand{\si}{{\sigma}}

\newcommand{\ov}{\overline}

\newcommand{\dd}{\mathrm{d}}

\DeclareMathOperator*{\esssup}{ess\,sup}

\makeatletter
\newcommand{\opnorm}{\@ifstar\@opnorms\@opnorm}
\newcommand{\@opnorms}[1]{%
	\left|\mkern-1.5mu\left|\mkern-1.5mu\left|
	#1
	\right|\mkern-1.5mu\right|\mkern-1.5mu\right|
}
\newcommand{\@opnorm}[2][]{%
	\mathopen{#1|\mkern-1.5mu#1|\mkern-1.5mu#1|}
	#2
	\mathclose{#1|\mkern-1.5mu#1|\mkern-1.5mu#1|}
}
\makeatother

\newtheoremstyle{neu}
{11pt}      
{11pt}      
{}                  
{}          
{\bfseries} 
{}          
{1em}  
{\textbf{\thmname{#1}\thmnumber{ #2}\thmnote{ (#3)}.}}          

\newtheoremstyle{proof}
{11pt}      
{11pt}      
{}                  
{}          
{\bfseries} 
{}            
{1em}          
{\textbf{\thmname{#1}\thmnote{ #3}.}}          

\newtheorem{Theorem}{Theorem}[section]

\newtheorem{Lemma}[Theorem]{Lemma}
\newtheorem{Proposition}[Theorem]{Proposition}
\theoremstyle{neu}

\newtheorem{Remark}[Theorem]{Remark}

\theoremstyle{proof}
\newtheorem{Proof}{Proof}

\numberwithin{equation}{section}

\allowdisplaybreaks[3]

\title{Normal approximation of the solution to the stochastic heat equation with Lévy noise}

\author{
	Carsten Chong\thanks{Institut de mathématiques, École Polytechnique Fédérale de Lausanne, Station 8, CH-1015 Lausanne, e-mail: carsten.chong@epfl.ch}
	~and
	Thomas Delerue\thanks{Chair of Mathematical Statistics, Technical University of Munich, Boltzmannstra\ss e 3, DE-85748 Garching, e-mail: thomas.delerue@tum.de}
}

\date{}

\begin{document}

\maketitle

\begin{abstract}
	Given a sequence $\dot{L}^\varepsilon$ of Lévy noises, we derive necessary and sufficient conditions in terms of their variances $\sigma^2(\varepsilon)$ such that the solution to the stochastic heat equation with noise $\sigma(\varepsilon)^{-1} \dot{L}^\varepsilon$ converges in law to the solution to the same equation with Gaussian noise. Our results apply to both equations with additive and multiplicative noise and hence lift the findings of S.\ Asmussen and J.\ Rosi\'nski [\textit{J. Appl. Probab.} \textbf{38} (2001) 482--493]  and S.\ Cohen and J.\ Rosi\'nski [\textit{Bernoulli} \textbf{13} (2007) 195--210] for finite-dimensional Lévy processes to the infinite-dimensional setting without making distributional assumptions on the solutions such as infinite divisibility. One important ingredient of our proof is to characterize the solution to the limit equation by a sequence of martingale problems. To this end, it is crucial to view the solution processes both as random fields and as càdlàg processes with values in a Sobolev space of negative real order.
\end{abstract}

\vfill

\noindent
\begin{tabbing}
	\> secondary: \, 46E35, 60G48
\end{tabbing}

\vspace{1cm}


\vspace{0.5cm}

\newpage

\section{Introduction}

The importance of the Gaussian distribution in probability theory and its popularity in applications are manifested in the central limit theorem: The total effect of a large number of small independent contributions is approximately normal. Therefore, when physical systems governed by one or several equations are pertubed by white noise (where ``white'' means stationary and uncorrelated), it is frequently assumed that the noise is \emph{Gaussian}. 

For example, in his Saint-Flour lecture notes \cite{Walsh}, J.\ B.\ Walsh discusses an application of parabolic stochastic PDEs to the modeling of neuron potentials. Subject to impulses arriving according to a marked Poisson point process (with mean $0$ and atoms of size $\dot L(t,x)$ at time $t$ and position $x$), the electrical potential $u(t,x)$ of the neuron, viewed as a thin cylinder of length, say, $\pi$, is then well described by the \emph{stochastic cable equation}
\begin{equation}\label{eq:cable-L} \partial_t u(t,x) = \partial_{xx} u(t,x) - u(t,x) + \dot L(t,x),\quad (t,x)\in[0,T]\times [0,\pi], \end{equation}
with suitable boundary and initial conditions. Arguing that ``the impulses are generally small, and there are many of them, so that in fact $\dot L$ is very nearly a white noise'' (\cite{Walsh}, p.\ 311), the author then approximates \eqref{eq:cable-L} by 
\begin{equation}\label{eq:cable-W}
\partial_t u(t,x) = \partial_{xx} u(t,x) - u(t,x) +\dot W(t,x),\quad (t,x)\in[0,T]\times [0,\pi],
\end{equation}
where $\dot W$ is a Gaussian space--time white noise.

But of course, the central limit theorem has limitations. In the absence of finite second moments, stable limits may arise; and if there are rare but large contributions, we may have a Poisson limit. In general, any \emph{infinitely divisible distribution} can arise as a possible limit of compound Poisson laws; see Corollary~3.8 in~\cite{Sato99}. This leads us to the following question: If we have a sequence of noises $\dot L^\varepsilon$ as above where the atom sizes of $\dot L^\varepsilon$ converge to $0$ as $\varepsilon\to0$, and if $u^\varepsilon$ denotes the solution to \eqref{eq:cable-L} with noise $\si(\varepsilon)^{-1}\dot L^\varepsilon$ (where $\si^2(\varepsilon)$ is the variance of $\dot L^\varepsilon$), do we have convergence in distribution of $u^\varepsilon$ to the solution $u$ of \eqref{eq:cable-W} with the Gaussian noise $\dot W$? A positive answer for this
\emph{normal approximation} is given in Theorem~7.10 in~\cite{Walsh}: If the atoms of $\dot L^\varepsilon$ are locally summable and the jump measure $Q^\varepsilon$ of $\dot L^\varepsilon$ satisfies
\begin{equation}\label{eq:deltacond} \frac{1}{\si^{2+\delta}(\varepsilon)}\int_\bbr |z|^{2+\delta}\,Q^\varepsilon(\dd z)\longrightarrow 0 \quad\text{as}\quad \varepsilon\to0\end{equation}
for some $\delta>0$, then $u^\varepsilon$ converges in distribution to $u$.

The purpose of this work is to substantially generalize this result in two aspects. Given that \eqref{eq:deltacond} is sufficient but not necessary for $u^\varepsilon \stackrel{d}{\longrightarrow} u$ (see Remark~\ref{rem:walsh} below), our first contribution is to show 
that the \emph{necessary and sufficient} condition for the normal approximation is  
\begin{equation} \label{AR condition}
\lim_{\varepsilon \rightarrow 0} 
\frac{1}{\sigma^2(\varepsilon)}  \int_{\vert z \vert > \kappa \sigma(\varepsilon)} z^2 \, Q^{\varepsilon}(\textrm{d}z) = 0
\end{equation}	
for all $\kappa>0$. In fact, if $L^\varepsilon$ (resp., $W$) is a Lévy process (resp., Brownian motion), the same condition was found to be necessary and sufficient for $\si(\varepsilon)^{-1} L^\varepsilon \stackrel{d}{\longrightarrow} W$ in \cite{CR}. Somewhat surprisingly, in the special case of \emph{small jump approximation}, that is, where $Q^\varepsilon(\textrm{d}z) = \mathbbm{1}_{\{|z|\leq\varepsilon\}}\,Q(\textrm{d}z)$ and $Q$ is a given Lévy measure, it was shown in \cite{AR} that condition \eqref{AR condition} \emph{fails} for prominent examples such as the compound Poisson or the gamma distribution. So in these cases, the small jump approximation is not true for Lévy processes and by our results, not true for stochastic PDEs, either.

Our second contribution is to consider equations with \emph{multiplicative} noise. To our best knowledge, previous works on the normal approximation of stochastic PDEs with jumps have only considered the situation of \emph{additive} noise; see, besides the mentioned results in \cite{Walsh}, also \cite{Kallianpur84, Walsh81} (there is, of course, literature concerning approximation of multiplicative Gaussian white noise by \emph{smoother} noises \cite{Bardina10,Hairer15}, but these problems are very different in nature than the one considered here). The proofs in \cite{Kallianpur84,Walsh81,Walsh} (as well as those of \cite{AR,CR}) are based on characteristic functions and the Lévy--Khintchine formula for infinitely divisible distributions, which obviously do not generalize to the situation of multiplicative noise.

Instead, our approach will be to show that $u^\varepsilon$ satisfies \emph{martingale problems} which, assuming \eqref{AR condition} only and not the stronger condition \eqref{eq:deltacond}, have a limit with a unique solution. But this leads to several complications. In order to prove  convergence of the associated martingales, we need some sort of uniformity in the time variable (as given, for example, by convergence in the Skorokhod topology). So taking simply the space $L^2([0,T]\times[0,\pi])$ to support the solutions $u^\varepsilon$ and $u$ will not be sufficient. This is why we will draw upon the results of \cite{CC} and view the solution $u^\varepsilon$ (and also $u$) as a càdlàg process on $[0,T]$ with values in the Sobolev space $H_{-r}$ for some $r>\frac12$ (see Section~\ref{section preliminaries} for a definition). In order to show tightness in that space with the Aldous criterion \cite{Aldous78}, we will use the factorization method from \cite{DaPrato87,Sanz} to obtain uniform bounds in time without taking moments of order higher than two. Another subtlety that arises in the analysis of the (semi-)martingales mentioned above is that their predictable characteristics are not given by a function of the former, which distinguishes our proof from the corresponding ones for (finite- or infinite-dimensional) stochastic differential equations in \cite{Fournier11,Kurtz91,Kurtz96}.

We shall also mention that the initial motivation in \cite{AR,CR} to study the normal approximation of Lévy processes comes from numerical simulation. Indeed, for stochastic PDEs as in \eqref{eq:cable-L} with multiplicative Lévy noise, the rate of convergence of a numerical scheme obtained by removing the small jumps of the noise is slower for noises with a high intensity of small jumps; see \cite{Chen16}. However, the results in~\cite{Kohatsu10} and~\cite{Kohatsu09} show that in the case of SDEs, an additional Gaussian approximation of the otherwise neglected small jumps improves the rate of convergence. We leave it to future research to examine to what extent this also holds for stochastic PDEs.

The remaining paper is organized as follows. In Section~\ref{sect:results}, we first describe in detail the considered equations and recall the definition of Sobolev spaces of real order in Section~\ref{section preliminaries} before we state our main result, Theorem~\ref{weak convergence}, in Section~\ref{section main result}. Here we also explain the main steps and ideas behind the proof, whereas the details are given in Section~\ref{sect:proof}.

\section{Results} \label{sect:results}

\subsection{Preliminaries} \label{section preliminaries}

Let $T > 0$ and consider on a filtered probability space $(\Omega, \mathcal{F}, \boldsymbol{F} = (\mathcal{F}_t)_{t \leq T}, \mathbb{P})$ that satisfies the usual conditions, for any  $\varepsilon > 0$, the \emph{stochastic heat equation} on $[0,T] \times [0,\pi]$ with Dirichlet boundary conditions:
\begin{equation} \label{SHE}
\left[
\begin{array}{ll}
\partial_t  u^{\varepsilon}(t,x) = \partial_{xx} u^{\varepsilon}(t,x) + f(u^{\varepsilon}(t,x))  
\frac{\displaystyle \dot{L}^{\varepsilon}(t,x)}{\displaystyle \sigma(\varepsilon)} , &
(t,x) \in [0, T] \times [0, \pi], \\
u^{\varepsilon}(t,0) = u^{\varepsilon}(t, \pi) = 0, & \textrm{for all} \,\, t \in [0,T], \\[1.6mm] \vspace{1.6mm}
u^{\varepsilon}(0,x) = 0, & \textrm{for all} \,\, x \in [0, \pi].
\end{array}	\right.
\end{equation}
The function $f \colon \mathbb{R} \longrightarrow \mathbb{R}$ in equation~\eqref{SHE} describes the \emph{multiplicative} part of the noise and will be assumed to be a Lipschitz continuous function. Concerning the driving noise $\sigma(\varepsilon)^{-1} \dot{L}^{\varepsilon}$, we assume that $\dot{L}^{\varepsilon}$ is a pure-jump Lévy space--time white noise on $[0,T] \times [0,\pi]$ given by
\begin{equation} \label{Levy noise}
	\begin{split}
		L^{\varepsilon}(A) = {}  &  
		\int_{0}^{T} \int_{0}^{\pi} \int_{\mathbb{R}} \mathbbm{1}_{A}(t,x)  z \, 
		(\mu^{\varepsilon} - \nu^{\varepsilon})(\textrm{d}t, \textrm{d}x, \textrm{d}z)
	\end{split}
\end{equation}
for all $A \in \mathcal{B}([0,T] \times [0,\pi])$. In this representation, $\mu^{\varepsilon}$ is a homogeneous Poisson random measure on $[0,T] \times [0,\pi] \times \mathbb{R}$ relative to the filtration $\boldsymbol{F}$, with intensity measure $\nu^{\varepsilon} = {\textrm{Leb}}_{[0,T] \times [0,\pi]} \otimes  Q^{\varepsilon}$. Here $Q^{\varepsilon}$ is a Lévy measure on $\mathbb{R}$, that is, $Q^{\varepsilon}(\{0\}) = 0$ and $\int_{\mathbb{R}} (1 \wedge z^2) \, Q^{\varepsilon}(\textrm{d}z) < \infty$. We refer to Chapter~II in \cite{JJ} for the definition of stochastic integrals with respect to Poisson random measures. Furthermore, we assume that for all $\varepsilon > 0$,
\begin{equation} \label{finite variance of Levy noise}
0<\sigma^2(\varepsilon) = \int_{\mathbb{R}} z^2 \, Q^{\varepsilon}(\textrm{d}z) < \infty.
\end{equation}
Note that this integral is the variance of $L^{\varepsilon}([0,1] \times [0,1])$. In the special case where we have a single Poisson random measure $\mu$ having intensity measure $\nu = {\textrm{Leb}}_{[0,T] \times [0,\pi]} \otimes Q$, setting
\begin{equation}\label{eq:Qvareps}
Q^{\varepsilon}(A) = \int_{\vert z \vert \leq \varepsilon} \mathbbm{1}_{A}(z) \, Q(\textrm{d}z), \quad A \in \mathcal{B}(\mathbb{R}), \quad \varepsilon > 0,
\end{equation}
leads us to the case of \emph{small jump approximation} considered in \cite{AR}.

A predictable random field $u^{\varepsilon} = \{u^{\varepsilon}(t,x) \mid (t,x) \in [0,T] \times [0,\pi] \}$ is called a \emph{mild solution} to~\eqref{SHE} if for all $(t,x) \in [0,T] \times [0,\pi]$,
\begin{equation} \label{random field solution Levy}
\begin{split}
u^{\varepsilon}(t,x) & = 
\int_0^t \int_{0}^{\pi} G_{t-s}(x,y) \frac{f(u^{\varepsilon}(s,y))}{\sigma(\varepsilon)} \, L^{\varepsilon}(\textrm{d}s, \textrm{d}y) \\ & =
\int_0^t \int_{0}^{\pi} \int_{\mathbb{R}} G_{t-s}(x,y)  f(u^{\varepsilon}(s,y)) \frac{z}{\sigma(\varepsilon)} \,
(\mu^{\varepsilon} - \nu^{\varepsilon})(\textrm{d}s, \textrm{d}y, \textrm{d}z)
\end{split}
\end{equation}
$\mathbb{P}$-almost surely, where
\begin{equation} \label{Green function}
G_t(x, y) = \frac{2}{\pi} \sum_{k=1}^{\infty}  \sin(kx) \sin(ky)  e^{-k^2t} \mathbbm{1}_{\{t \geq 0\}},
\end{equation}
for $(t,x,y) \in [0,T] \times [0,\pi]^2$, is the Dirichlet Green's function of the heat operator on $[0,\pi]$. 

The existence of a mild solution $u^{\varepsilon}$ to~\eqref{SHE} is guaranteed by Theorem~3.1 in~\cite{CC2} and condition \eqref{finite variance of Levy noise} on the Lévy measure $Q^{\eps}$ and it is, up to modifications, unique among all predictable random fields satisfying
\begin{equation} \label{bound second moment of solution}
\sup_{(t,x) \in [0,T] \times [0,\pi]} 
\mathbb{E}\left[ {\vert u^{\varepsilon}(t,x) \vert}^p \right] < \infty
\end{equation}
for any $0 < p \leq 2$ and $\eps > 0$.

In this paper, we want to examine when the \emph{normal approximation} holds for $u^{\varepsilon}$, that is, when $u^{\varepsilon}$ can be approximated in law by the mild solution $u$ to the same stochastic heat equation as above, but driven by a Gaussian space--time white noise on $[0,T] \times [0,\pi]$:
\begin{equation} \label{SHE 2} 
\left[
\begin{array}{ll}
\partial_t u(t,x) = \partial_{xx} u(t,x) + f(u(t,x))  \dot{W}(t,x), &
(t,x) \in [0, T] \times [0, \pi], \\
u(t,0) = u(t, \pi) = 0, & \textrm{for all} \,\, t \in [0,T], \\
u(0,x) = 0, & \textrm{for all} \,\, x \in [0, \pi].
\end{array}	\right.
\end{equation}
The driving noise $\dot{W}$ is now a centered Gaussian random field $\left \{ W(A) \mid A \in \mathcal{B}([0,T] \times [0,\pi]) \right \}$ with covariance structure $\mathbb{E}[W(A)W(B)]={\textrm{Leb}}_{[0,T] \times [0,\pi]}(A \cap B)$ for any measurable sets $A, B \subseteq [0,T] \times [0,\pi]$. It is well-known (see, for example, Theorem~3.2 in~\cite{Walsh}) that, up to modifications, equation~\eqref{SHE 2} has a unique mild solution $u$ satisfying the corresponding bound in~\eqref{bound second moment of solution} for all $p > 0$.

Throughout this work, we will look at the mild solutions $u^{\varepsilon}$ and $u$ from two different points of view. First, they are random elements in the function space $L^2([0,T] \times [0,\pi])$ as the uniform bound~\eqref{bound second moment of solution} shows. But then, as mentioned in the introduction, we will need stronger path regularity in the time variable for our proofs. This is why we shall consider $u^{\varepsilon}$ and $u$ also as stochastic processes with values in an infinite dimensional space, which we will describe in the following. 

Consider for any $r > 0$, the fractional Sobolev space
\begin{equation*}
H_{r}([0,\pi]) = \left \{ \phi \in L^2([0,\pi]) \, \, \Big \vert \, \, 
\sum_{k=1}^{\infty} (1 + k^2)^{r} \langle \phi, \phi_k \rangle^2 < \infty \right \},
\end{equation*}
where $\phi_k(x) = \sqrt{2/\pi} \sin(kx)$, $k \in \mathbb{N}$, form an orthonormal basis of $L^2([0,\pi])$. This is a Hilbert space with scalar product
\begin{equation*}
{\langle f, g \rangle}_r = 
\sum_{k=1}^{\infty} (1 + k^2)^{r}  
\langle f, \phi_k \rangle \langle g, \phi_k \rangle, \quad f, g \in H_{r}([0,\pi]),
\end{equation*}
and norm ${\Vert \phi \Vert}_r = \sqrt{\langle \phi, \phi \rangle_r}$ for $\phi \in H_{r}([0,\pi])$. 

The topological dual $H_{-r}([0,\pi])$ of $H_{r}([0,\pi])$ is also a Hilbert space, whose dual norm ${\Vert \cdot \Vert}_{-r}$ can be expressed, by the Riesz representation theorem, as
\begin{equation} \label{dual norm}
{\Vert \phi' \Vert}_{-r}^2 = \sum_{k=1}^{\infty} 
\langle \phi', \phi_{r,k} \rangle^2 = \sum_{k=1}^{\infty} (1 + k^2)^{-r} \langle \phi', \phi_{k} \rangle^2, \quad \phi' \in H_{-r}([0, \pi]).
\end{equation}
(Note that if $\varphi_1, \varphi_2$ are elements of the same $L^2$-space, then $\langle \varphi_1, \varphi_2 \rangle$ will always denote the standard scalar product of that space. If $\phi$ is an element of a Hilbert space and $\phi'$ an element of its topological dual, then $\langle \phi', \phi \rangle$ will always denote the dual pairing of $\phi'$ with $\phi$.)

Coming back to the mild solution to~\eqref{SHE}, if we identify $u^{\varepsilon}$ with the process ${(u^{\varepsilon}_t)}_{t \leq T}$ where
\begin{equation} \label{identification of solution}
u^{\varepsilon}_t \colon H_{r}([0,\pi])  \longrightarrow \mathbb{R}, \quad
\phi \mapsto \langle u^{\varepsilon}( t, \cdot), \phi \rangle 
= \int_{0}^{\pi} u^{\varepsilon}( t, y) \phi(y) \, \textrm{d}y
\end{equation}
for all $t \leq T$, then by Theorem 2.5 in \cite{CC}, $u^{\varepsilon}$ has a càdlàg modification in $H_{-r}([0,\pi])$ for any $r > 1/2$, which will be denoted by $\overline{u}^{\varepsilon} = (\overline{u}^{\varepsilon}_t)_{t \leq T}$ throughout this work. Similarly, by the identification~\eqref{identification of solution} and Corollary~3.4 in \cite{Walsh}, the mild solution $u$ to~\eqref{SHE 2} has a continuous modification $\overline{u}$ in $H_{-r}([0,\pi])$ for each $r > 1/2$.

\subsection{Main result} \label{section main result}

We now introduce the Cartesian product
\begin{equation} \label{cartesian product}
\Omega^* = L^2([0,T] \times [0,\pi]) \times D([0,T], H_{-r}([0,\pi])),  
\end{equation}
with $r > 1/2$. Let $d_1$ denote the metric induced by the $L^2$-norm on $L^2([0,T] \times [0,\pi])$ and $d_2$ be the Skorokhod metric on $D([0,T], H_{-r}([0,\pi]))$. We then equip $\Omega^*$ with the product metric
\begin{equation} \label{product metric}
\tau((f_1, x_1), (f_2, x_2)) = d_1(f_1, f_2) + d_2(x_1, x_2)
\end{equation}
for any $f_1, f_2 \in L^2([0,T] \times [0,\pi])$ and $x_1, x_2 \in D([0,T], H_{-r}([0,\pi]))$. The main result of this paper is the following limit theorem.

\begin{Theorem} \label{weak convergence}
	Assume that $L^\varepsilon$ is given by \eqref{Levy noise} with a variance $\sigma^2(\varepsilon)$ that satisfies \eqref{finite variance of Levy noise} for all $\varepsilon >0$. 
Let $u^{\varepsilon}$ be the $L^2([0,T] \times [0,\pi])$-valued mild solution to the stochastic heat equation~\eqref{SHE} driven by $\dot{L}^{\varepsilon}/\sigma(\varepsilon)$ and $\overline{u}^{\varepsilon}$ be its càdlàg modification in $H_{-r}([0,\pi])$. 	
Similarly, let $u$ be the $L^2([0,T] \times [0,\pi])$-valued mild solution to the stochastic heat equation~\eqref{SHE 2} driven by $\dot{W}$ and $\overline{u}$ be its continuous modification in $H_{-r}([0,\pi])$. 

Suppose also that the Lipschitz function $f$ in \eqref{SHE} satisfies $f(0) \neq 0$. Then, as $\varepsilon \rightarrow 0$,
\begin{equation} \label{goal}
(u^{\varepsilon}, \overline{u}^{\varepsilon}) \stackrel{d}{\longrightarrow} (u, \overline{u}) \quad \textrm{in} \quad
(\Omega^*, \tau)
\end{equation}
for all $r > 1/2$ if and only if~\eqref{AR condition} holds for all $\kappa > 0$.
\end{Theorem}

\begin{Remark}
	We can generalize Theorem~\ref{weak convergence} to nonzero initial conditions. Assume that in both equations~\eqref{SHE} and~\eqref{SHE 2} we now have $u^\varepsilon(0,x)=u(0,x) = u_0(x)$ for all $x \in [0, \pi]$, where $u_0 \colon [0,\pi] \longrightarrow \mathbb{R}$ is a bounded continuous function with $u_0(0) = u_0(\pi) = 0$. Define	
	\begin{equation*}
	u_0(t,x) = \int_{0}^{\pi} G_t(x,y) u_0(y) \, \textrm{d}y, \quad (t,x) \in [0, T] \times [0, \pi].
	\end{equation*}	
	Then Theorem~\ref{weak convergence} can be shown in a completely analogous manner if we assume that there exists $(t_0, x_0) \in [0, T] \times [0, \pi]$ such that $f(u_0(t_0,x_0)) \neq 0$ (instead of $f(0) \neq 0$), and this assumption is only needed for showing the necessity of~\eqref{AR condition}. To be more precise, since the mild solution to~\eqref{SHE 2} now satisfies	
	\begin{equation*}
	u(t,x) = u_0(t,x) +
	\int_0^t \int_{0}^{\pi} G_{t-s}(x,y) f(u(s,y)) \, W(\textrm{d}s, \textrm{d}y) 
	\end{equation*}
	$\mathbb{P}$-almost surely, a similar argument as in Remark~\ref{Remark on f} shows that $\mathbb{P}(f(u(t_1,x_1)) \neq 0) > 0$ for some $(t_1, x_1) \in [0, T] \times [0, \pi]$ and hence, the expectation in~\eqref{NR5} is nonzero.	
\end{Remark}

\begin{Remark}\label{rem:walsh} Let us relate the two conditions \eqref{eq:deltacond} and \eqref{AR condition} to each other. Using Hölder's and Chebyshev's inequalities, we see from the estimate
	\begin{equation*}
	\begin{split}
		\frac{1}{\sigma^2(\varepsilon)}  \int_{\vert z \vert > \kappa \sigma(\varepsilon)} z^2 \, Q^{\varepsilon}(\textrm{d}z)
		&\leq 	\frac{1}{\sigma^2(\varepsilon)}  \left(\int_\mathbb{R} |z|^{2+\delta} \, Q^{\varepsilon}(\textrm{d}z)\right)^{\frac{2}{2+\delta}} Q^{\varepsilon}(\{\vert z \vert > \kappa \sigma(\varepsilon)\})^{\frac{\delta}{2+\delta}}\\
		&\leq \frac{1}{\kappa^{\delta}\sigma^{2+\delta}(\varepsilon)}  \int_\mathbb{R} |z|^{2+\delta} \, Q^{\varepsilon}(\textrm{d}z)
	\end{split}
	\end{equation*}
	that \eqref{eq:deltacond} implies \eqref{AR condition}.
	
	The other implication is not true in general. For example, assume that $Q^\varepsilon$ has density
		\begin{equation*} q_\varepsilon(z)=\frac{1}{2\vert z \vert^2}\mathbbm{1}_{\{ \vert z \vert\leq\varepsilon\}} + \frac{\varepsilon^2}{2C\vert z \vert^3\log(1+\vert z \vert)^2}\mathbbm{1}_{\{ \vert z \vert>1\}},\quad   z  \in \mathbb{R}, \end{equation*}
	where $C=\int_1^\infty  z^{-1}\log(1+ z)^{-2}\,\textrm{d}z$.
	Then $\int_\mathbb{R} \vert z \vert^{2+\delta} \,Q^\varepsilon(\textrm{d}z)=\infty$ for every $\varepsilon,\delta>0$, so condition \eqref{eq:deltacond} does not hold. But a direct calculation shows that $\sigma^2(\varepsilon)=\varepsilon + \varepsilon^2$. So given $\kappa>0$,  we have $1>\kappa\sigma(\varepsilon)>\kappa\sqrt{\varepsilon}>\varepsilon$ for small values of $\varepsilon$, which implies \eqref{AR condition} because
\begin{equation*} \lim_{\varepsilon \rightarrow 0} 
\frac{1}{\sigma^2(\varepsilon)}  \int_{\vert z \vert > \kappa \sigma(\varepsilon)} z^2 \, Q^{\varepsilon}(\textrm{d}z) = \lim_{\varepsilon \rightarrow 0} 
\frac{1}{\varepsilon+\varepsilon^2}  \int_{\vert z \vert > \varepsilon} z^2 \, Q^{\varepsilon}(\textrm{d}z) =\lim_{\varepsilon \rightarrow 0} 
\frac{\varepsilon^2}{\varepsilon+\varepsilon^2} =0.  \end{equation*}
\end{Remark}

\begin{Proof}[of Theorem~\ref{weak convergence}]
We begin by showing that~\eqref{AR condition} implies the weak convergence~\eqref{goal}. Since $\Omega^*$ is a metric space, we follow the classical scheme of first showing tightness and then uniqueness of the limiting distribution. 

In Theorem~\ref{tightness in L2}, we show that $\{ u^{\varepsilon} \mid \varepsilon > 0 \}$ is tight in $L^2([0,T] \times [0,\pi])$ and in Theorem~\ref{tightness in Skorokhod space} that $\{ \overline{u}^{\varepsilon} \mid \varepsilon > 0 \}$ is tight in $D([0,T], H_{-r}([0,\pi]))$. By the subsequence principle, this immediately implies that the random elements $\{ (u^{\varepsilon}, \overline{u}^{\varepsilon}) \mid \varepsilon >0 \}$ are tight in $(\Omega^*, \tau)$. As it turns out, no assumptions on the Lévy noise $\dot{L}^{\varepsilon}$ other than the ones specified in~\eqref{Levy noise} and~\eqref{finite variance of Levy noise} are needed for this tightness property.

As a consequence, we can apply Prokhorov's theorem, which provides for any sequence ${(\varepsilon_k)}_{k \in \mathbb{N}}$ with $\varepsilon_k \rightarrow 0$, a subsequence ${(\varepsilon_{k_l})}_{l \in \mathbb{N}}$ such that $(u^{\varepsilon_{k_l}}, \overline{u}^{\varepsilon_{k_l}})_{l \in \mathbb{N}}$ converges weakly to some distribution on $(\Omega^*, \tau)$ as $l \rightarrow \infty$. For notational simplicity, we will assume without loss of generality that the whole sequence $(u^{\varepsilon_k}, \overline{u}^{\varepsilon_k})_{k \in \mathbb{N}}$ converges weakly. 

Since $(\Omega^*, \tau)$ is a complete separable metric space, we can further apply Skorokhod's representation theorem (see Theorem~4.30 in~\cite{Kallenberg}) and obtain random elements 
\begin{equation} \label{Skorokhod representation 1}
(v^{k}, \overline{v}^{k}), (v, \overline{v})  \colon (\overline{\Omega}, \overline{\mathcal{F}}, \overline{\mathbb{P}}) \longrightarrow (\Omega^*, \tau),
\end{equation}
defined on a possibly different probability space $(\overline{\Omega}, \overline{\mathcal{F}}, \overline{\mathbb{P}})$,
satisfying the following properties:
\begin{equation} \label{Skorokhod representation 3}
\begin{split}
& (v^{k}, \overline{v}^{k}) \stackrel{d}{=} (u^{\varepsilon_{k}}, \overline{u}^{\varepsilon_{k}}) \quad \textrm{for all} \quad k \in \mathbb{N} \quad \textrm{and} \\ & 
(v^{k}, \overline{v}^{k})(\overline{\omega}) \longrightarrow (v, \overline{v})(\overline{\omega}) \quad \textrm{in} \quad (\Omega^*, \tau) \quad \textrm{as} \quad k \rightarrow \infty \quad \textrm{for all} \quad \overline{\omega} \in \overline{\Omega}.
\end{split}
\end{equation}

We will show that 
\begin{equation*}
(v, \overline{v}) \stackrel{d}{=} (u, \overline{u}),
\end{equation*}
which in turn implies~\eqref{goal}. To do this, we first define a filtration $\overline{\boldsymbol{F}} = (\overline{\mathcal{F}}_t)_{t \leq T}$ on $\overline{\Omega}$ by setting
\begin{equation} \label{def filtration Om Stern}
\overline{\mathcal{F}}_t = \bigcap_{u > t} \sigma \left(v^k(s,x), \overline{v}^k_s \mid 0 \leq x \leq \pi,\, k \in \mathbb{N},\, s \leq u \right), \quad t \leq T. 
\end{equation}
We further define for $\xi \in \mathbb{R}$, $\phi \in C_c^{\infty}((0,\pi))$ and $t \leq T$,
\begin{equation} \label{characteristics of eta}
\begin{split}
& \overline{B}_t = \int_{0}^{t} \langle v(s, \cdot), \phi'' \rangle \, \textrm{d}s, \quad
\overline{C}_t = 
\int_{0}^{t} \int_{0}^{\pi} f^2(v(s,x))  \phi^2(x) \, \textrm{d}s \, \textrm{d}x,  \quad
\overline{A}_t = i \xi \overline{B}_t - \frac{1}{2}  \xi^2 \overline{C}_t
\end{split}
\end{equation}
as well as
\begin{equation} \label{def ov M}
\overline{M}_t = e^{i  \xi \langle \overline{v}_t, \phi \rangle}  - 
\int_{0}^{t} e^{i  \xi \langle \overline{v}_s, \phi \rangle} \, \overline{A} (\textrm{d}s).
\end{equation}
We will then show in Theorem~\ref{convergence of the MPs} that, under assumption~\eqref{AR condition}, the pair $(v, \overline{v})$ satisfies the following martingale problem. For all $\xi \in \mathbb{R}$ and $\phi \in C_c^{\infty}((0,\pi))$, the process $(\overline{M}_t)_{t \leq T}$ is a martingale with respect to $(\overline{\Omega}, \overline{\mathcal{F}}, \overline{\boldsymbol{F}}, \overline{\mathbb{P}})$. Note that we are able to obtain this property only because in~\eqref{Skorokhod representation 3}, we have $\overline{\omega}$-wise convergence both in the Skorokhod topology and in $L^2([0,T] \times [0,\pi])$, which is the reason why we view the solutions to~\eqref{SHE} and~\eqref{SHE 2} as \emph{pairs} in $\Omega^*$.


Next, we will show in Theorem~\ref{solution of MP is solution of SHE with W} that this martingale property in turn implies that there exists a Gaussian space--time white noise $\dot{\widetilde{W}}$ on $[0,T] \times [0,\pi]$, possibly defined on a filtered extension $(\widetilde{\Omega}, \widetilde{\mathcal{F}}, \widetilde{\boldsymbol{F}}, \widetilde{\mathbb{P}})$ of $(\overline{\Omega}, \overline{\mathcal{F}}, \overline{\boldsymbol{F}}, \overline{\mathbb{P}})$ such that, with probability one, the random field $v$ is equal in $L^2([0,T] \times [0,\pi])$ to the mild solution $\widetilde{v}$ to the stochastic heat equation
\begin{equation} \label{SHE 2'} 
\left[
\begin{array}{ll}
\partial_t  \widetilde{v}(t,x) = \partial_{xx} \widetilde{v}(t,x) + f(\widetilde{v}(t,x))  \dot{\widetilde{W}}(t,x), &
(t,x) \in [0, T] \times [0, \pi], \\
\widetilde{v}(t,0) = \widetilde{v}(t, \pi) = 0, & \textrm{for all} \,\, t \in [0,T], \\
\widetilde{v}(0,x) = 0, & \textrm{for all} \,\, x \in [0, \pi],
\end{array}	\right.
\end{equation}
and such that $\overline{v}$ is indistinguishable from the continuous version in $H_{-r}([0,\pi])$ of $\widetilde{v}$, which concludes the first part of the proof.

For the second part, under the assumption $f(0) \neq 0$, Theorem~\ref{necessary condition} directly shows that~\eqref{goal} implies~\eqref{AR condition}.
\qed
\end{Proof}

\section{Details of the proof} \label{sect:proof}

In the remainder of this work, the letter $C$ will always denote a strictly positive constant whose value may change from line to line. Furthermore, by the Lipschitz continuity of the function $f$, there exists a positive constant $K$ that we hold fixed from now on such that $\vert f(x) \vert \leq K \vert x \vert + \vert f(0) \vert $ for all $x \in \mathbb{R}$.

\subsection{Tightness} \label{sect:tight}

We start with three lemmas that will provide uniform bounds in $\varepsilon > 0$ for the second moments of $u^{\varepsilon}$, which will be crucial for proving tightness of $\{ (u^{\varepsilon}, \overline{u}^{\varepsilon}) \mid \varepsilon > 0 \}$ in $(\Omega^*, \tau)$.

\begin{Lemma} \label{Lemma uniform bound}
	The family $\{ u^{\varepsilon} \mid \varepsilon >0 \}$ of mild solutions to (\ref{SHE}) satisfies 
	\begin{equation} \label{uniform bound second moment}
	\sup_{\varepsilon > 0}  \sup_{(t,x) \in [0,T] \times [0,\pi]} 
	\mathbb{E}\left[ \vert u^{\varepsilon}(t,x) \vert^2 \right] < \infty
	\end{equation}
	and this uniform bound only depends on the Lipschitz function $f$.
\end{Lemma}

\begin{Proof}
Using It\={o}'s isometry and the definition~\eqref{random field solution Levy} of $u^{\varepsilon}$, we have for fixed $\varepsilon > 0$ and $(t,x) \in [0,T] \times [0,\pi]$,
\begin{equation*}
\begin{split}
\mathbb{E} \left[ {\vert u^{\varepsilon}(t,x) \vert}^2 \right] & = 
\mathbb{E} \left[
\int_0^t \int_{0}^{\pi} \int_{\mathbb{R}} G^2_{t-s}(x,y) \frac{f^2(u^{\varepsilon}(s,y))}{\sigma^2(\varepsilon)}  z^2 \, 
\nu^{\varepsilon}(\textrm{d}s, \textrm{d}y, \textrm{d}z) \right] \\ & =
\mathbb{E} \left[ \int_0^t \int_{0}^{\pi} G_{t-s}^2(x,y)
f^2(u^{\varepsilon}(s,y)) \, \textrm{d}s \, \textrm{d}y \right]
\left(\frac{1}{\sigma^2(\varepsilon)} 
\int_{\mathbb{R}} z^2 \, Q^{\varepsilon}(\textrm{d}z)\right)  \\ & = 
\int_0^t \int_{0}^{\pi} G_{t-s}^2(x,y) 
\mathbb{E} \left[ f^2(u^{\varepsilon}(s,y)) \right]  \textrm{d}s \, \textrm{d}y.
\end{split}
\end{equation*}
Using the Lipschitz continuity of $f$ and the elementary inequality $(a + b)^2 \leq 2a^2 + 2b^2$, we then obtain for $\varepsilon > 0$ and $(t,x) \in [0,T] \times [0,\pi]$,
\begin{equation} \label{NR7}
\mathbb{E} \left[ {\vert u^{\varepsilon}(t,x) \vert}^2 \right] \leq
C \int_0^t \int_{0}^{\pi}  G^2_{t-s}(x,y)  \mathbb{E} \left[ \vert u^{\varepsilon}(s,y) \vert^2 \right] \textrm{d}s \, \textrm{d}y +
C \int_0^t \int_{0}^{\pi}  G^2_{t-s}(x,y) \, \textrm{d}s \, \textrm{d}y.
\end{equation}
Now in order to find a bound for $\mathbb{E} [ {\vert u^{\varepsilon}(t,x) \vert}^2 ]$ uniformly in $t$, $x$ and $\varepsilon$, we will use a comparison principle for deterministic Volterra equations. By (B.5) in \cite{Bally}, there exists a constant $C > 0$ such that $\vert G_t(x,y) \vert \leq C g_t(x-y)$ on $[0,T] \times [0,\pi]^2$, where
\begin{equation} \label{Green's function on R}
g_t(x) = \frac{1}{\sqrt{4\pi t}} \exp \left(- \frac{\vert x \vert^2}{4 t} \right) \mathbbm{1}_{\{ t \geq 0 \}}
\end{equation}
is the heat kernel on $\mathbb{R}$. Since $\int_{0}^{T} \int_{\mathbb{R}} {\vert g_t(x) \vert}^q \, \textrm{d}t \, \textrm{d}x < \infty$ for all $q < 3$, we obtain
\begin{equation*}
\sup_{(t,x) \in [0,T] \times [0,\pi]} \int_0^t \int_{0}^{\pi} G^2_{t-s}(x,y)  \, \textrm{d}s \, \textrm{d}y \leq 
C \int_{0}^{T} \int_{\mathbb{R}} g^2_t(x)  \, \textrm{d}t \, \textrm{d}x < \infty.
\end{equation*}
Recall from (\ref{bound second moment of solution}) that $\mathbb{E} [ {\vert u^{\varepsilon}(t,x) \vert}^2 ]$ is uniformly bounded in $(t,x)$ for fixed $\varepsilon > 0$. Therefore, by Lemma 6.4 (2) and (3) in \cite{CC2}, the mild solution $u^{\varepsilon}$ satisfies
\begin{equation*}
\mathbb{E} \left[ {\vert u^{\varepsilon}(t,x) \vert}^2 \right] \leq v(t,x)
\end{equation*}
for all $(t,x) \in [0,T] \times [0,\pi]$ and $\varepsilon > 0$, where $v$ is the unique nonnegative solution of the deterministic Volterra equation
\begin{equation*}
\begin{split}
v(t,x) & = 
C \int_0^t \int_{0}^{\pi} G^2_{t-s}(x,y) v(s,y) \, \textrm{d}s \, \textrm{d}y \\ & \quad \,\, +
C \int_0^t \int_{0}^{\pi} G^2_{t-s}(x,y) \, \textrm{d}s \, \textrm{d}y, \quad 
(t,x) \in [0,T] \times [0,\pi],
\end{split}
\end{equation*}
and satisfies $\sup_{(t,x) \in [0,T] \times [0,\pi]}  v(t,x) < \infty$. \qed
\end{Proof}

The next lemma gives an alternative integral representation of $u^{\varepsilon}$ and is an extension of the factorization method in~\cite{DaPrato87} and~\cite{Sanz}.

\begin{Lemma} \label{Lemma factorization method}
	For $\delta \in (0, 1/4)$ define 
	\begin{equation*}
	Y_{\delta}^{\varepsilon}(t,x) = \int_{0}^{t} \int_{0}^{\pi}  \frac{G_{t-s}(x,y)}{(t-s)^{\delta}}  \frac{f(u^{\varepsilon}(s,y))}{\sigma(\varepsilon)} \, L^{\varepsilon}(\mathrm{d}s, \mathrm{d}y), \quad (t,x) \in [0,T] \times [0,\pi].
	\end{equation*}
	We then have 
	\begin{equation} \label{uniform bound second moment 2}
	\sup_{\varepsilon > 0} \sup_{(t,x) \in [0,T] \times [0,\pi]} 
	\mathbb{E}\left[ \vert Y_{\delta}^{\varepsilon}(t,x) \vert^2 \right] < \infty
	\end{equation}
	and for all $(t,x) \in [0,T] \times [0,\pi]$, the representation
	\begin{equation} \label{factorization method}
	u^{\varepsilon}(t,x) = \frac{\sin(\delta \pi)}{\pi} \int_{0}^{t} \int_{0}^{\pi} \frac{G_{t-s}(x,y)}{(t-s)^{1-\delta}}  Y_{\delta}^{\varepsilon}(s,y) \, \mathrm{d}s \, \mathrm{d}y
	\end{equation}
	holds $\mathbb{P}$-almost surely.
\end{Lemma}

\begin{Proof}
First, using It\={o}'s isometry, the Lipschitz continuity of $f$ and~Lemma \ref{Lemma uniform bound}, we have
\begin{equation*}
\begin{split}
\mathbb{E} \left[ \vert Y_{\delta}^{\varepsilon}(t,x) \vert^2 \right] & = 
\mathbb{E} \left[
\int_0^t \int_{0}^{\pi} \frac{G^2_{t-s}(x,y)}{(t-s)^{2\delta}} 
f^2(u^{\varepsilon}(s,y)) \, \textrm{d}s \, \textrm{d}y \right] \left(\frac{1}{\sigma^2(\varepsilon)}
\int_{\mathbb{R}} z^2 \, Q^{\varepsilon}(\textrm{d}z)\right) \\ & =
\int_{0}^{t} \int_{0}^{\pi} \frac{G^2_{t-s}(x,y)}{(t-s)^{2\delta}} \mathbb{E} \left[ f^2(u^{\varepsilon}(s,y)) \right] \,\textrm{d}s \, \textrm{d}y \leq 
C \int_{0}^{t} \int_{0}^{\pi} \frac{G^2_{t-s}(x,y)}{(t-s)^{2\delta}} \, \textrm{d}s \, \textrm{d}y.
\end{split}
\end{equation*}
The last integral on the right-hand side is finite if $\delta < 1/4$. Indeed, by (B.5) in \cite{Bally}, it can be bounded by
\begin{equation*}
C \int_{0}^{t} \int_{0}^{\pi} \frac{1}{(t-s)^{2\delta  + 1}} \exp \left(- \frac{\vert x - y \vert^2}{t-s} \right) \, \textrm{d}s \, \textrm{d}y = 
C \int_{0}^{t} \frac{1}{(t-s)^{2\delta  + 1/2}} \, \textrm{d}s = Ct^{-2\delta + 1/2}.
\end{equation*}
The identity~\eqref{factorization method} follows in the same way as Lemma 5 in \cite{Sanz}. \qed 
\end{Proof}

\begin{Lemma} \label{Lemma uniform bound moment of norm}
	The family $\{ u^{\varepsilon} \mid \varepsilon >0 \}$ of mild solutions to (\ref{SHE}) satisfies 
	\begin{equation} \label{uniform bound moment of norm}
	\sup_{\varepsilon > 0}  \mathbb{E}\left[ 
	\left( \int_{0}^{T} \left(\int_{0}^{\pi} \vert u^{\varepsilon}(t,x) \vert^2 \, \mathrm{d}x \right)^p \mathrm{d}t \right)^{1/p}
	\right] < \infty
	\end{equation}
	for all $p \in (1, 4/3)$.
\end{Lemma}

\begin{Proof}
We will use the integral representation~\eqref{factorization method} of Lemma \ref{Lemma factorization method}. Fix $\delta \in (0, 1/4)$ and $t \in [0,T]$. Using Fubini's theorem, we have
\begin{equation*}
\begin{split}
& \int_{0}^{\pi} \vert u^{\varepsilon}(t,x) \vert^2 \, \textrm{d}x \\ & \qquad  =
C  \int_{0}^{\pi}  \left( 
\int_{0}^{t} \int_{0}^{\pi}  \int_{0}^{t} \int_{0}^{\pi} 
\frac{G_{t-s}(x,y)}{(t-s)^{1-\delta}}   
\frac{G_{t-s'}(x,y')}{(t-s')^{1-\delta}} Y_{\delta}^{\varepsilon}(s,y) Y_{\delta}^{\varepsilon}(s',y') \, \textrm{d}y' \, \textrm{d}s' \, \textrm{d}y \, \textrm{d}s 
\right)\, \textrm{d}x.
\end{split}
\end{equation*}
By the semigroup property of the Green's function, the integral on the right-hand side is equal to
\begin{equation*}
\begin{split}
& \int_{0}^{t} \int_{0}^{\pi} \int_{0}^{t} \int_{0}^{\pi} 
\frac{G_{2t-s-s'}(y,y')}{(t-s)^{1-\delta} (t-s')^{1-\delta}} Y_{\delta}^{\varepsilon}(s,y)  Y_{\delta}^{\varepsilon}(s',y')  \, \textrm{d}y' \, \textrm{d}s' \, \textrm{d}y \, \textrm{d}s \\ & \qquad = 
2 \int_{0}^{t} \int_{0}^{\pi} \int_{0}^{s} \int_{0}^{\pi} 
\frac{G_{2t-s-s'}(y,y')}{(t-s)^{1-\delta} (t-s')^{1-\delta}} Y_{\delta}^{\varepsilon}(s,y)  Y_{\delta}^{\varepsilon}(s',y')  \, \textrm{d}y' \, \textrm{d}s' \, \textrm{d}y \, \textrm{d}s.
\end{split}
\end{equation*}
Using (B.5) in \cite{Bally} and \eqref{Green's function on R}, we obtain
\begin{equation*}
\int_{0}^{\pi} \vert u^{\varepsilon}(t,x) \vert^2 \, \textrm{d}x \leq 
C \int_{0}^{t} \int_{0}^{\pi} \int_{0}^{s} \int_{0}^{\pi} 
\frac{g_{2t-s-s'}(y-y')}{(t-s)^{1-\delta} (t-s')^{1-\delta}} \vert Y_{\delta}^{\varepsilon}(s,y)  Y_{\delta}^{\varepsilon}(s',y') \vert \, \textrm{d}y' \, \textrm{d}s' \, \textrm{d}y \, \textrm{d}s.
\end{equation*}
Now let $p \in (1, 4/3)$, take the $\Vert \cdot \Vert_{L^p([0,T])}$-norm of $t \mapsto \int_{0}^{\pi} \vert u^{\varepsilon}(t,x) \vert^2 \, \textrm{d}x$ and apply Minkowski's integral inequality to obtain
\begin{equation*}
\begin{split}
& \left( \int_{0}^{T} \left(\int_{0}^{\pi} \vert u^{\varepsilon}(t,x) \vert^2 \, \textrm{d}x \right)^p\, \textrm{d}t \,\right)^{1/p} \\ & \qquad \leq C
\int_{0}^{T}  \int_{0}^{\pi} \int_{0}^{s} \int_{0}^{\pi} 
\left(\int_{0}^{T} 
\left(\frac{g_{2t-s-s'}(y-y') \mathbbm{1}_{\{ s \leq t \}}}{(t-s)^{1-\delta} \, (t-s')^{1-\delta}}\right)^p \, \textrm{d}t\right)^{1/p} 
\vert Y_{\delta}^{\varepsilon}(s,y)  Y_{\delta}^{\varepsilon}(s',y') \vert \, \textrm{d}y' \, \textrm{d}s' \, \textrm{d}y \, \textrm{d}s.
\end{split}
\end{equation*}
Take expectation, use the Cauchy--Schwarz inequality and \eqref{uniform bound second moment 2} to further obtain
\begin{equation} \label{NR12}
\begin{split}
& \mathbb{E} \left[\left( \int_{0}^{T} \left(\int_{0}^{\pi} \vert u^{\varepsilon}(t,x) \vert^2 \, \textrm{d}x \right)^p\, \textrm{d}t \right)^{1/p}\right] \\ & \qquad
 \leq C
\int_{0}^{T} \int_{0}^{\pi} \int_{0}^{s} \int_{0}^{\pi} 
\left(\int_{s}^{T} 
\left(\frac{g_{2t-s-s'}(y-y')}{(t-s)^{1-\delta}  (t-s')^{1-\delta}}\right)^p  \textrm{d}t \right)^{1/p} \, \textrm{d}y' \, \textrm{d}s' \, \textrm{d}y \, \textrm{d}s.
\end{split}
\end{equation}
Note that the right-hand side of \eqref{NR12} does not depend on $\varepsilon$ anymore. We now consider the integrand
\begin{equation*}
\int_{s}^{T} 
\left(\frac{g_{2t-s-s'}(y-y')}{(t-s)^{1-\delta}  (t-s')^{1-\delta}}\right)^p  \,\textrm{d}t = 
\int_{0}^{T-s} 
\frac{g^{p}_{2t + s-s'}(y-y')}{t^{(1-\delta)p}  (t+s-s')^{(1-\delta)p}} \, \textrm{d}t.
\end{equation*}
For fixed $x \in \mathbb{R}$, the maximum of the function $t \mapsto g_t(x)$ is $C/\vert x \vert$ for some $C$ that is independent of $x$. Let $\eta \in (0,1)$ and consider the estimate 
\begin{equation*}
g_t(x) = g_t(x)^{1 - \eta}  g_t(x)^{\eta} \leq C \frac{1}{\vert x \vert^{1 - \eta}}  \frac{1}{t^{\eta/2}}, \quad t > 0, \quad x \in \mathbb{R}.
\end{equation*}
Since $s' \leq s$, we obtain
\begin{equation*}
\begin{split}
& \left(\int_{0}^{T-s} 
\frac{g^{p}_{2t + s-s'}(y-y')}{t^{(1-\delta)p}  (t+s-s')^{(1-\delta)p}} \, \textrm{d}t \right)^{1/p} 
\\ & \qquad   \leq C  \frac{1}{\vert y -y ' \vert^{1 - \eta}}
\left( \int_{0}^{T} \frac{1}{t^{(1-\delta)p}  (t+s-s')^{(1-\delta)p}  (2t + s-s')^{p\eta/2}} \, \textrm{d}t \right)^{1/p} 
\\ & \qquad 
\leq C \frac{1}{\vert y -y ' \vert^{1 - \eta}}
\left( \int_{0}^{T} \frac{1}{t^{((1-\delta)+\eta/2)p}  (t+s-s')^{(1-\delta)p}} \, \textrm{d}t \right)^{1/p}.
\end{split}
\end{equation*}
Moreover, the integral $\int_{0}^{\pi} \int_{0}^{\pi} \vert y -y ' \vert^{\eta - 1} \, \textrm{d}y' \, \textrm{d}y$ is finite because $\eta > 0$, and the expectation in \eqref{NR12} is now bounded by 
\begin{equation} \label{NR13}
C \int_{0}^{T} \int_{0}^{s} \left(\int_{0}^{T} 
\frac{1}{t^{((1-\delta)+\eta/2)p}  (t+s-s')^{(1-\delta)p}} \, \textrm{d}t \right)^{1/p} \,\textrm{d}s' \, \textrm{d}s
\end{equation}
for any $\delta \in (0, 1/4)$ and $\eta \in (0,1)$. By assumption, $3/4 < 1/p < 1$ and $3/4 < (1 - \delta) + \eta/2 < 3/2$. Hence, we can choose $\delta$ and $\eta$ such that $(1 - \delta) + \eta/2 < 1/p$. As a consequence, by Lemma 2 of Chapter 1 in \cite{Friedman}, the estimate
\begin{equation*}
\int_{0}^{T} 
\frac{1}{t^{((1-\delta)+\eta/2)p} (t+s-s')^{(1-\delta)p}} \, \textrm{d}t \leq C
(s - s')^{1 - ((1-\delta)+\eta/2)p - (1-\delta)p}
\end{equation*}
holds for $0 \leq s' < s \leq T$. We have assumed that $((1-\delta)+\eta/2)p + (1-\delta)p > 1$ since otherwise, the last integral is bounded by some constant, which immediately implies that the expectation in \eqref{NR12} is uniformly bounded in $\varepsilon$. 

Therefore, we further estimate the integral in \eqref{NR13} by
\begin{equation*}
C \int_{0}^{T} \int_{0}^{s} (s - s')^{1/p - 2(1-\delta)-\eta/2} \, \textrm{d}s' \, \textrm{d}s = 
C \int_{0}^{T} \int_{0}^{s} r^{1/p - 2(1-\delta)-\eta/2} \, \textrm{d}r \, \textrm{d}s, 
\end{equation*}
which is finite because our choice of $\delta$ and $\eta$ implies $2(1-\delta) + \eta/2 - 1/p < 1$. This concludes the proof. \qed 
\end{Proof}

We can now proceed to showing tightness.

\begin{Theorem} \label{tightness in L2}
	The family $\{ u^{\varepsilon} \mid \varepsilon >0 \}$ of mild solutions to (\ref{SHE}) is tight in the Hilbert space $L^2([0,T] \times [0,\pi])$.
\end{Theorem}

\begin{Proof}
	It is easy to see that the functions  
	\begin{equation*}
	\psi_{ij}(t,x) = \overline{\phi}_i(t) \phi_j(x), \quad (t,x) \in [0,T] \times [0,\pi],
	\end{equation*}
	where $\overline{\phi}_i(t) = \sqrt{2/T}  \sin (it  \pi /T )$ and $\phi_j(x) = \sqrt{2/\pi} \sin(jx)$ for all $i, j \in \mathbb{N}$, form an orthonormal basis of $L^2([0,T] \times [0,\pi])$.	
	
	First, using the stochastic Fubini theorem (see, for example, Theorem 2.6 in \cite{Walsh}), we have for all $i,j \in \mathbb{N}$,
	\begin{equation} \label{NR9}
	\begin{split}
	\langle u^{\varepsilon}, \psi_{ij} \rangle & = 
	\int_{0}^{T} \int_{0}^{\pi} \left(\int_0^t \int_{0}^{\pi} G_{t-s}(x,y) \frac{f(u^{\varepsilon}(s,y))}{\sigma(\varepsilon)} \, L^{\varepsilon}(\textrm{d}s, \textrm{d}y) \right) \psi_{ij}(t,x) \, \textrm{d}t \, \textrm{d}x \\ & =
	\frac{1}{\sigma(\varepsilon)} \int_0^T \int_{0}^{\pi} f(u^{\varepsilon}(s,y)) 
	\left(\int_{s}^{T} \int_{0}^{\pi} G_{t-s}(x,y) \psi_{ij}(t,x) \, \textrm{d}t \, \textrm{d}x \right)  \, L^{\varepsilon}(\textrm{d}s, \textrm{d}y).
	\end{split}
	\end{equation}
	Define for all $i,j \in \mathbb{N}$,
	\begin{equation*}
	H_{ij}(s,y) = \int_{s}^{T} \int_{0}^{\pi} G_{t-s}(x,y)  \psi_{ij}(t,x) \, \textrm{d}t \, \textrm{d}x, \quad (s,y) \in [0,T] \times [0,\pi].
	\end{equation*}
	Using Fubini's theorem, the expression~\eqref{Green function} of the Green's function $G$ and the orthogonal properties of $\phi_j$, we obtain for all $(s,y) \in [0,T] \times [0,\pi]$,
	\begin{equation} \label{NR23}
	H_{ij}(s,y) =
	\int_{s}^{T} \overline{\phi}_i(t) \left(\int_{0}^{\pi} G_{t-s}(x,y) \phi_j(x) \, \textrm{d}x \right) \,\textrm{d}t =
	\int_{s}^{T} \overline{\phi}_i(t) \phi_j(y) e^{-j^2(t-s)} \, \textrm{d}t.
	\end{equation}	
	Using the integral formula 
	$\int (\sin ax)  e^{bx} \, \textrm{d}x =  \left(b \sin ax - a \cos ax \right) e^{bx}/ (a^2 + b^2) + C$, we can further calculate
	\begin{equation*}
	\begin{split}
	H_{ij}(s,y) & = 
	\sqrt{\frac{2}{T}}  \phi_j(y)  e^{j^2 s} 
	\int_{s}^{T} \sin \left(i \frac{\pi}{T}  t\right) e^{-j^2 t} \, \textrm{d}t \\ & =
	\sqrt{\frac{2}{T}}  \phi_j(y) 
	\frac{1}{i^2 (\pi/T)^2 + j^4} 
	\left( e^{-j^2 (T-s)} i  \frac{\pi}{T}  (-1)^{i+1} + j^2 \sin \left(i \frac{\pi}{T} s\right) +
	i \frac{\pi}{T} \cos \left(i \frac{\pi}{T} s  \right) \right) \\ & \leq
	C\left( \frac{i}{i^2 + j^4} + \frac{j^2}{i^2 + j^4} \right) \leq C  \frac{1}{i + j^2} \quad \textrm{for all} \quad i,j \in \mathbb{N}.
	\end{split}
	\end{equation*}
	For the $L^2$-norm of $H_{ij}$, we then have
	\begin{equation*}
	\int_{0}^{T} \int_{0}^{\pi} H_{ij}^2(s,y) \, \textrm{d}s \, \textrm{d}y \leq C  \frac{1}{i^2 + j^4}
	\end{equation*}
	for all $i,j \in \mathbb{N}$. Since
	\begin{equation*}
	\sum_{i,j =2}^{\infty} \frac{1}{i^2 + j^4} \leq \int_{1}^{\infty} \int_{1}^{\infty} \frac{1}{x^2 + y^4} \, \mathrm{d}x \, \mathrm{d}y = 
	\int_{1}^{\infty} \frac{\arctan y^2}{y^2} \, \mathrm{d}y < \infty,
	\end{equation*}
	we obtain from Lemma~\ref{Lemma uniform bound}
	\begin{equation*}
	\sum_{i, j=1}^{\infty} \sup_{\varepsilon > 0} \mathbb{E} \left[ \langle u^{\varepsilon}, \psi_{ij} \rangle^2 \right] \leq 
	C \sum_{i, j=1}^{\infty} \int_{0}^{T} \int_{0}^{\pi} H_{ij}^2(s,y) \, \textrm{d}s \, \textrm{d}y < \infty,
	\end{equation*}
	which implies that
	\begin{equation*}
	\begin{split}
	& \sup_{\varepsilon > 0}  \mathbb{P} \left( \sum_{i, j \geq N} \langle u^{\varepsilon}, \psi_{ij} \rangle^2 > \delta \right) \leq 
	\frac{1}{\delta} \sum_{i, j \geq N}^{\infty} \sup_{\varepsilon > 0} \mathbb{E} \left[ \langle u^{\varepsilon}, \psi_{ij} \rangle^2 \right] \longrightarrow 0 \quad \textrm{as} \quad N \rightarrow \infty
	\end{split}
	\end{equation*}
	for all $\delta > 0$. Moreover, again by Lemma~\ref{Lemma uniform bound}, we have
	\begin{equation*}
	\begin{split}
	\sup_{\varepsilon > 0}  \mathbb{P} \left( \sum_{i, j < N} \langle u^{\varepsilon}, \psi_{ij} \rangle^2 > \delta \right) & \leq 
	\frac{1}{\delta}  \sup_{\varepsilon > 0}  \mathbb{E} \left[ \sum_{i, j < N} \langle u^{\varepsilon}, \psi_{ij} \rangle^2 \right] \leq 
	\frac{1}{\delta}  \sup_{\varepsilon > 0} \mathbb{E} \left[ \sum_{i, j =1}^{\infty} \langle u^{\varepsilon}, \psi_{ij} \rangle^2 \right] \\ & =
	\frac{1}{\delta}  \sup_{\varepsilon > 0}  \mathbb{E} \left[ 
	\int_{0}^{T} \int_{0}^{\pi} u^{\varepsilon}(t,x)^2 \, \textrm{d}t \, \textrm{d}x  \right] \leq \frac{C}{\delta} 
	\longrightarrow 0
	\end{split}
	\end{equation*}
	as $\delta \rightarrow \infty$ for all $N \in \mathbb{N}$. Therefore, we can conclude from Theorem 1 in \cite{Suquet} that $\{ u^{\varepsilon} \mid \varepsilon >0 \}$ is tight in $L^2([0,T] \times [0,\pi])$. \qed
\end{Proof}

The next two propositions will imply that $\{ \overline{u}^{\varepsilon} \mid \varepsilon >0 \}$ is tight in $D([0,T], H_{-r}([0,\pi]))$. 

\begin{Proposition} \label{Aldous condition}
	The càdlàg processes $\{ \overline{u}^{\varepsilon} \mid \varepsilon >0 \}$ satisfy the Aldous condition: Let $(\varepsilon_n)_{n \in \mathbb{N}}$ and $(h_n)_{n \in \mathbb{N}}$ be sequences of positive numbers with $\varepsilon_n \rightarrow 0$ and $h_n \rightarrow 0$ as $n \rightarrow \infty$. For each $n \in \mathbb{N}$, further let $\tau_n  \in [0,T]$ be a stopping time with respect to the filtration generated by the stochastic process $(\overline{u}_t^{\varepsilon_n})_{t \leq T}$. Then we have for any $r> 1/2$,
	\begin{equation*}
	\mathbb{E}\left[
	\Vert \overline{u}^{\varepsilon_n}_{\tau_n + h_n} - \overline{u}^{\varepsilon_n}_{\tau_n} \Vert_{-r}^2
	\right] \longrightarrow 0 \quad \textrm{as} \quad  n \rightarrow \infty.
	\end{equation*}
\end{Proposition}

\begin{Proof}
Recall the expression of the dual norm ${\Vert \cdot \Vert}_{-r}$ in \eqref{dual norm}. We have
\begin{equation} \label{NR8}
\Vert \overline{u}^{\varepsilon_n}_{\tau_n + h_n} - \overline{u}^{\varepsilon_n}_{\tau_n} \Vert_{-r}^2  =
\sum_{k=1}^{\infty} ( 1 + k^2 )^{-r} 
\left(\langle \overline{u}^{\varepsilon_n}_{\tau_n + h_n}, \phi_{k} \rangle - \langle \overline{u}^{\varepsilon_n}_{\tau_n }, \phi_{k} \rangle\right)^2.
\end{equation}
We will find a convenient semimartingale decomposition for the real-valued stochastic process $\langle \overline{u}^{\varepsilon}, \phi_{k} \rangle = \left(\langle \overline{u}^{\varepsilon}_{t}, \phi_{k} \rangle\right)_{t \leq T}$ for any $\varepsilon > 0$ and $k \in \mathbb{N}$ that will then allow us to estimate the expectation of the terms appearing in \eqref{NR8}. 

First, proceeding as in~\eqref{NR9} and~\eqref{NR23}, we have for all $t \leq T$,
\begin{equation*}
\begin{split}
\langle u^{\varepsilon}(t, \cdot), \phi_k \rangle & = 
\int_{0}^{\pi} \left(\int_0^t \int_{0}^{\pi} G_{t-s}(x,y) \frac{f(u^{\varepsilon}(s,y))}{\sigma(\varepsilon)} \, L^{\varepsilon}(\textrm{d}s, \textrm{d}y) \right) \phi_k(x) \, \textrm{d}x \\ & =
\frac{1}{\sigma(\varepsilon)} \int_0^t \int_{0}^{\pi} f(u^{\varepsilon}(s,y)) \phi_k(y) e^{-k^2(t-s)} \,  L^{\varepsilon}(\textrm{d}s, \textrm{d}y).
\end{split}
\end{equation*}
If we define 
\begin{equation*}
X^{k, \varepsilon}_t = \int_{0}^{t} \int_{0}^{\pi} f(u^{\varepsilon}(s,y)) \phi_k(y) \, L^{\varepsilon}(\textrm{d}s, \textrm{d}y), \quad t \leq T,
\end{equation*}
for all $k \in \mathbb{N}$ and $\varepsilon > 0$, then
\begin{equation*}
\frac{1}{\sigma(\varepsilon)} \int_0^t \int_{0}^{\pi} f(u^{\varepsilon}(s,y)) \phi_k(y) e^{-k^2(t-s)} \,
L^{\varepsilon}(\textrm{d}s, \textrm{d}y) = \frac{1}{\sigma(\varepsilon)} 
\int_{0}^{t} e^{-k^2(t-s)}  \,X^{k, \varepsilon}(\textrm{d}s)
\end{equation*}
for all $t \leq T$, where the last term is the It\={o} integral of the deterministic function $s \mapsto e^{-k^2(t-s)}$ against the square-integrale martingale $X^{k, \varepsilon}$. Because the integrand is a $C^{\infty}$-function, the integration by parts formula for semimartingales yields
\begin{equation*}
\int_{0}^{t} e^{-k^2(t-s)} \,X^{k, \varepsilon}(\textrm{d}s) = 
X^{k, \varepsilon}_t - 
\int_{0}^{t} X^{k, \varepsilon}_s k^2 e^{-k^2(t-s)} \, \textrm{d}s.
\end{equation*}
Altogether we obtain the semimartingale decomposition
\begin{equation} \label{semimartingale decomposition}
\begin{split}
\langle u^{\varepsilon}(t, \cdot), \phi_k \rangle & =
\int_{0}^{t} \int_{0}^{\pi} \frac{f(u^{\varepsilon}(s,y))}{\sigma(\varepsilon)} \phi_k(y) \, L^{\varepsilon}(\textrm{d}s, \textrm{d}y) \\ & \quad \,\, -   
\int_{0}^{t} \left( \int_{0}^{s} \int_{0}^{\pi} \frac{f(u^{\varepsilon}(r,y))}{\sigma(\varepsilon)} \phi_k(y) \, L^{\varepsilon}(\textrm{d}r, \textrm{d}y) \right) k^2 e^{-k^2(t-s)} \, \textrm{d}s
\end{split}
\end{equation}
for all $t \leq T$, $k \in \mathbb{N}$ and $\varepsilon > 0$.

Now the process $\langle \overline{u}^{\varepsilon}, \phi_k \rangle$ is the càdlàg version of $(\langle u^{\varepsilon}(t, \cdot), \phi_k \rangle)_{t \leq T}$, so we can infer that $\langle \overline{u}^{\varepsilon}, \phi_k \rangle$ and the right-hand side of \eqref{semimartingale decomposition} are indistinguishable since the latter is also càdlàg. Coming back to \eqref{NR8}, we can now decompose
\begin{equation*}
\langle \overline{u}^{\varepsilon_n}_{\tau_n + h_n}, \phi_k \rangle - \langle \overline{u}^{\varepsilon_n}_{\tau_n}, \phi_k \rangle = I_{k,n} + J_{k,n}^1 + J_{k,n}^2,
\end{equation*}
where
\begin{equation*}
\begin{split}
I_{k,n} & = \int_0^{T} \int_{0}^{\pi} \frac{f(u^{\varepsilon_n}(s,y))}{\sigma(\varepsilon_n)} \phi_k(y)  
\mathbbm{1}_{(\tau_n, \tau_n + h_n]}(s) \,
L^{\varepsilon_n}(\textrm{d}s, \textrm{d}y), \\
J_{k,n}^1 & = 
\int_{0}^{\tau_n} \left(
\int_0^s \int_{0}^{\pi} \frac{f(u^{\varepsilon_n}(r,y))}{\sigma(\varepsilon_n)} \phi_k(y) \,
L^{\varepsilon_n}(\textrm{d}r, \textrm{d}y)
\right) k^2 \left(e^{-k^2(\tau_n-s)} - e^{-k^2(\tau_n + h_n -s)}\right)\, \textrm{d}s, \\
J_{k,n}^2 & = -
\int_{\tau_n}^{\tau_n + h_n} \left(
\int_0^s \int_{0}^{\pi} \frac{f(u^{\varepsilon_n}(r,y))}{\sigma(\varepsilon_n)} \phi_k(y) \,
L^{\varepsilon_n}(\textrm{d}r, \textrm{d}y)
\right) k^2 e^{-k^2(\tau_n + h_n -s)} \, \textrm{d}s
\end{split}
\end{equation*}
for all $k,n \in \mathbb{N}$. We now gather some moment estimates for these three terms. First, for the martingale term $I_{k,n}$, we have by It\={o}'s isometry,
\begin{equation} \label{NR6}
\begin{split}
\mathbb{E}[I^2_{k,n}] & = 
\mathbb{E} \left[ 
\int_0^{T} \int_{0}^{\pi} \int_{\mathbb{R}} 
f^2(u^{\varepsilon_n}(s,y)) \phi_k^2(y) 
\mathbbm{1}_{(\tau_n, \tau_n + h_n]}(s) 
\frac{z^2}{\sigma^2(\varepsilon_n)} \, \textrm{d}s \, \textrm{d}y \, Q^{\varepsilon_n}(\textrm{d}z) 
\right] \\ & =
\mathbb{E} \left[ 
\int_0^{T} \int_{0}^{\pi} 
f^2(u^{\varepsilon_n}(s,y)) \phi_k^2(y) 
\mathbbm{1}_{(\tau_n, \tau_n + h_n]}(s) \, \textrm{d}s \, \textrm{d}y 
\right].
\end{split}
\end{equation}
Using the Lipschitz continuity of $f$, we can bound the last term in \eqref{NR6} by
\begin{equation*}
C \mathbb{E} \left[ 
\int_0^{T} \int_{0}^{\pi} 
u^{\varepsilon_n}(s,y)^2  
\mathbbm{1}_{(\tau_n, \tau_n + h_n]}(s) \, \textrm{d}s \, \textrm{d}y 
\right] + 
C \mathbb{E} \left[ \int_0^{T} \int_{0}^{\pi} \mathbbm{1}_{(\tau_n, \tau_n + h_n]}(s) \, \textrm{d}s \, \textrm{d}y \right]
\end{equation*}
for any $k \in \mathbb{N}$. The second term equals $C \pi h_n$, which converges to 0 as $n \rightarrow \infty$. 

For the first term, choose $p \in (1, 4/3)$. Using Hölder's inequality and Lemma~\ref{Lemma uniform bound moment of norm}, we obtain
\begin{equation*}
\begin{split}
& \mathbb{E} \left[ 
\int_0^{T} \left(\int_{0}^{\pi} 
u^{\varepsilon_n}(s,y)^2 \, \textrm{d}y \right)
\mathbbm{1}_{(\tau_n, \tau_n + h_n]}(s) \, \textrm{d}s 
\right] \\ & \qquad  \leq 
\mathbb{E} \left[
\left(\int_0^{T} \left(\int_{0}^{\pi} 
u^{\varepsilon_n}(s,y)^2 \, \textrm{d}y \right)^p\,  \textrm{d}s \right)^{1/p}  
\left( \int_0^{T} \mathbbm{1}_{(\tau_n, \tau_n + h_n]}(s) \, \textrm{d}s  \right)^{1-1/p}
\right] \\ & \qquad \leq h_n^{1-1/p} \sup_{\varepsilon > 0}  
\mathbb{E} \left[
\left(\int_0^{T} \left(\int_{0}^{\pi} 
u^{\varepsilon}(s,y)^2 \, \textrm{d}y \right)^p \, \textrm{d}s \right)^{1/p} \right] \longrightarrow 0  \quad \textrm{as} \quad n \rightarrow \infty.
\end{split}
\end{equation*}
Altogether, this implies $\mathbb{E} [ I_{k,n}^2 ] \rightarrow 0$ as $n \rightarrow \infty$ for all $k \in \mathbb{N}$.

Next, we have
\begin{equation*}
\begin{split}
\int_{0}^{\tau_n} k^2 \left(e^{-k^2(\tau_n-s)} - e^{-k^2(\tau_n + h_n -s)}\right) \, \textrm{d}s  & = 
\left(1 - e^{-k^2  h_n} \right) \int_{0}^{\tau_n} k^2  e^{-k^2(\tau_n-s)} \, \textrm{d}s \\ & =
\left(1 - e^{-k^2  h_n} \right) \left( 1 - e^{-k^2  \tau_n} \right) \leq 1 - e^{-k^2 h_n},
\end{split}
\end{equation*}
and by It\={o}'s isometry as well as Lemma~\ref{Lemma uniform bound},
\begin{equation*}
\mathbb{E} \left[ 
\left( 
\int_0^{T} \int_{0}^{\pi} \frac{f(u^{\varepsilon_n}(s,y))}{\sigma(\varepsilon_n)}  \phi_k(y) \, 
L^{\varepsilon_n}(\textrm{d}s, \textrm{d}y)
\right)^2
\right] =
\int_0^{T} \int_{0}^{\pi} \mathbb{E} \left[ 
f^2(u^{\varepsilon_n}(s,y)) \right] \phi_k^2(y) \,  \textrm{d}s \, \textrm{d}y  \leq C
\end{equation*}
for all $k, n \in \mathbb{N}$. Therefore, by Doob's inequality, $\mathbb{E} [ ( J_{k,n}^1 )^2 ]$ is bounded by
\begin{equation*}
\begin{split}
& \mathbb{E} \Bigg[
\sup_{s \leq T} \Bigg \vert 
\int_0^s \int_{0}^{\pi} \frac{f(u^{\varepsilon_n}(r,y))}{\sigma(\varepsilon_n)}  \phi_k(y) \,
L^{\varepsilon_n}(\textrm{d}r, \textrm{d}y) \Bigg \vert^2  \Bigg(
\int_{0}^{\tau_n}  k^2 \left(e^{-k^2(\tau_n-s)} - e^{-k^2(\tau_n + h_n -s)}\right) \, \textrm{d}s  \Bigg)^2 \Bigg] \\ & \qquad \leq
\left(1 - e^{-k^2 h_n} \right)^2 
\mathbb{E} \Bigg[
\sup_{s \leq T} \Bigg \vert 
\int_0^s \int_{0}^{\pi} \frac{f(u^{\varepsilon_n}(r,y))}{\sigma(\varepsilon_n)}  \phi_k(y) \,
L^{\varepsilon_n}(\textrm{d}r, \textrm{d}y) \Bigg \vert^2 \Bigg] \\ & \qquad \leq
C \left(1 - e^{-k^2  h_n} \right)^2 \longrightarrow 0 \quad \textrm{as} \quad n \rightarrow \infty
\end{split}
\end{equation*}
for all $k \in \mathbb{N}$. 

Finally, with $\int_{\tau_n}^{\tau_n + h_n} k^2 e^{-k^2(\tau_n + h_n -s)} \, \textrm{d}s = 1 - e^{-k^2  h_n}$ and similar calculations, $\mathbb{E} [ ( J_{k,n}^2)^2 ]$ is bounded by
\begin{equation*}
\begin{split}
& \mathbb{E} \Bigg[
\sup_{s \leq T} \Bigg \vert 
\int_0^s \int_{0}^{\pi} \frac{f(u^{\varepsilon_n}(r,y))}{\sigma(\varepsilon_n)}  \phi_k(y) \,
L^{\varepsilon_n}(\textrm{d}r, \textrm{d}y) \Bigg \vert^2 \Bigg(
\int_{\tau_n}^{\tau_n + h_n} k^2 e^{-k^2(\tau_n + h_n -s)} \, \textrm{d}s \Bigg)^2
\Bigg] \\ & \qquad = \left( 1 - e^{-k^2  h_n} \right)^2 
\mathbb{E} \Bigg[
\sup_{s \leq T} \Bigg \vert 
\int_0^s \int_{0}^{\pi} \frac{f(u^{\varepsilon_n}(r,y))}{\sigma(\varepsilon_n)}  \phi_k(y) \,
L^{\varepsilon_n}(\textrm{d}r, \textrm{d}y) \Bigg \vert^2 \Bigg] \\ & \qquad \leq
C \left( 1 - e^{-k^2  h_n} \right)^2 \longrightarrow 0 \quad \textrm{as} \quad n \rightarrow \infty
\end{split}
\end{equation*}
for all $k \in \mathbb{N}$. As a consequence, recalling that $r > 1/2$ and thus  $\sum_{k=1}^{\infty} {(1 + k^2)}^{-r} < \infty$, we obtain by \eqref{NR8} and dominated convergence,
\begin{equation*}
\begin{split}
\mathbb{E} \left[
\Vert \overline{u}^{\varepsilon_n}_{\tau_n + h_n} - \overline{u}^{\varepsilon_n}_{\tau_n} \Vert_{-r}^2 \right] & =
\sum_{k=1}^{\infty} ( 1 + k^2 )^{-r}  \mathbb{E} \left[
\left(\langle \overline{u}^{\varepsilon_n}_{\tau_n + h_n}, \phi_{k} \rangle - \langle \overline{u}^{\varepsilon_n}_{\tau_n }, \phi_{k} \rangle\right)^2 \right] \\ & \leq
3 \sum_{k=1}^{\infty} ( 1 + k^2 )^{-r} \left(\mathbb{E} \left[ I_{k,n}^2 \right] + \mathbb{E} \left[ ( J_{k,n}^1)^2 \right] + 
\mathbb{E} \left[ ( J_{k,n}^2 )^2 \right] \right)  \longrightarrow 0 
\end{split}
\end{equation*}
as $n \rightarrow \infty$, which is the assertion of the proposition. \qed
\end{Proof}

\begin{Proposition} \label{Second condition for tightness}
	For any fixed $t \leq T$ and $r > 1/2$, the random elements $\{ \overline{u}_t^{\varepsilon} \mid \varepsilon > 0 \}$ are tight in $H_{-r}([0,\pi])$.
\end{Proposition}

\begin{Proof} Proceeding as in the proof of Proposition \ref{Aldous condition} and using Lemma~\ref{Lemma uniform bound}, we have
\begin{equation*}
\mathbb{E} \left[ {\langle \overline{u}^{\varepsilon}_t, \phi_k \rangle}^2 \right] =
\int_0^t \int_{0}^{\pi} \mathbb{E} \left[ 
f^2(u^{\varepsilon}(s,y)) \right]
\phi^2_k(y)  e^{-2k^2(t-s)} \, \textrm{d}s \, \textrm{d}y \leq  C \frac{1}{k^2} \left( 1 - e^{-2k^2 t} \right)
\end{equation*}
for any $k \in \mathbb{N}$ and $t \leq T$. Hence, we have for all $q < 1/2$, $t \leq T$ and $\varepsilon > 0$,
\begin{equation*}
\begin{split}
& \mathbb{E} \left[ {\Vert \overline{u}_t^{\varepsilon} \Vert}^2_{q} \right] = 
\sum_{k=1}^{\infty} (1 + k^2)^{q} \mathbb{E} \left[ {\langle \overline{u}^{\varepsilon}_t, \phi_k \rangle}^2 \right] \leq 
C \sum_{k=1}^{\infty} (1 + k^2)^{q} \frac{1}{k^2} < \infty,
\end{split}
\end{equation*}
and thus $\overline{u}_t^{\varepsilon} \in H_q([0,\pi])$ $\mathbb{P}$-almost surely. 

Because the penultimate term in the inequality above does not depend on $\varepsilon$, by Markov's inequality, we can further deduce 
\begin{equation} \label{NR24}
\lim_{\delta \rightarrow \infty}  \sup_{\varepsilon > 0}  \mathbb{P} \left( {\Vert \overline{u}_t^{\varepsilon} \Vert}_{q} > \delta \right) = 0
\end{equation}
for all $q < 1/2$ and $t \leq T$. Since the embeddings 
\begin{equation*}
H_q([0,\pi]) \hookrightarrow L^2([0,\pi]) \hookrightarrow H_{-r}([0,\pi])
\end{equation*}
are compact for $0 < q < 1/2 < r$ by Theorem 4.58 in \cite{Demengel}, it follows that $\{ \ov{u}_t^{\varepsilon} \mid \varepsilon >0 \}$ is tight in $H_{-r}([0,\pi])$ for any fixed $t \leq T$ and $r > 1/2$. \qed 
\end{Proof}

\begin{Theorem} \label{tightness in Skorokhod space}
	 For any $r> 1/2$, the càdlàg modifications $\{ \overline{u}^{\varepsilon} \mid \varepsilon >0 \}$ are tight in the Skorokhod space $D([0,T], H_{-r}([0,\pi]))$.
\end{Theorem}

\begin{Proof}
By Theorem 6.8 in \cite{Walsh}, this is a direct consequence of Propositions~\ref{Aldous condition}~and~\ref{Second condition for tightness}. \qed

\end{Proof}

\subsection{Characterization of the limit} \label{sect:charlim}

After proving tightness in Section~\ref{sect:tight}, our next goal is to characterize the limit distribution of weakly converging subsequences. Following the outline of the proof of Theorem~\ref{weak convergence}, the first step is to show that under condition~\eqref{AR condition} on the Lévy measure $Q^{\varepsilon}$, the process $\overline{M}$ in~\eqref{def ov M} is a martingale with respect to the filtration $\overline{\boldsymbol{F}}$ defined in~\eqref{def filtration Om Stern}. In order to achieve this result, which is Theorem~\ref{convergence of the MPs} below, we prove that the pairs $(u^{\varepsilon}, \overline{u}^{\varepsilon})$ satisfy related martingale problems (Theorem~\ref{solution of SHE with Lévy noise is solution of MP}) and that these ``converge'' as $\varepsilon \rightarrow 0$ (Theorem~\ref{convergence ov M k to ov M}).

Recall that for all test functions $\phi \in C_c^{\infty}((0, \pi))$ and fixed $t \leq T$,
\begin{equation} \label{weak formulation}
\begin{split}
& \int_{0}^{\pi} u^{\varepsilon}(t,x)  \phi(x) \, \textrm{d}x = 
\int_{0}^{t} \int_{0}^{\pi} u^{\varepsilon}(s,x)  \phi''(x) \, \textrm{d}s \, \textrm{d}x + 
\int_{0}^{t} \int_{0}^{\pi} 
\frac{f(u^{\varepsilon}(s,x))}{\sigma(\varepsilon)}  \phi(x) \,
L^{\varepsilon}(\textrm{d}s, \textrm{d}x)
\end{split}
\end{equation}
$\mathbb{P}$-almost surely. This follows, in a similar way to Theorem 3.2 in \cite{Walsh}, from the fact that in our situation, we may apply the stochastic Fubini theorem; see, for example, Theorem 2.6 in \cite{Walsh}.

\begin{Theorem} \label{solution of SHE with Lévy noise is solution of MP}
	For each $\varepsilon > 0$, the pair $(u^{\varepsilon}, \overline{u}^{\varepsilon})$ where $u^{\varepsilon} \in L^2([0,T] \times [0,\pi])$ is the mild solution to the stochastic heat equation (\ref{SHE}) and $\overline{u}^{\varepsilon}$ is its càdlàg modification in $H_{-r}([0,\pi])$, with $r > 1/2$, satisfies the following martingale problem. For all $\xi \in \mathbb{R}$ and $\phi \in C_c^{\infty}((0,\pi))$, the complex-valued stochastic process
	\begin{equation} \label{definition M eps}
	\begin{split}
	M_t^{\varepsilon} & = e^{i  \xi \langle \overline{u}^{\varepsilon}_t,  \phi \rangle}  -
	i \xi \int_{0}^{t} e^{i  \xi \langle \overline{u}^{\varepsilon}_s, \phi \rangle}  \langle u^{\varepsilon}(s, \cdot), \phi'' \rangle \, \mathrm{d}s \\ & \quad \,\, - 
	\int_{0}^{t} \int_{0}^{\pi} \int_{\mathbb{R}} 
	e^{i \xi \langle \overline{u}^{\varepsilon}_{s}, \phi \rangle} 
	\left(e^{i  \xi  \frac{f(u^{\varepsilon}(s,x))}{\sigma(\varepsilon)} \phi(x)  z} - 1 -
	i\xi \frac{f(u^{\varepsilon}(s,x))}{\sigma(\varepsilon)} \phi(x)  z \right)\,
	 \mathrm{d}s \, \mathrm{d}x \,  Q^{\varepsilon}(\mathrm{d}z), 
	\end{split}
	\end{equation}
	with $t \leq T$, is a square-integrable $\boldsymbol{F}$-martingale with the uniform bound
	\begin{equation} \label{bound second moment M eps}
	\sup_{\varepsilon > 0}  \sup_{t \leq T}  \mathbb{E} \left[ \big \vert M^{\varepsilon}_t \big \vert^2 \right] < \infty.
	\end{equation}
\end{Theorem}

\begin{Proof}
	First, since $\overline{u}^{\varepsilon}$ is the càdlàg version of $u^{\varepsilon}$, the stochastic process $\langle \overline{u}^{\varepsilon}, \phi \rangle$ is indistinguishable from the right-hand side of~\eqref{weak formulation} for any $\phi \in C_c^{\infty}((0,\pi))$. This directly implies 
	that $\langle \overline{u}^{\varepsilon}, \phi \rangle$ is an $\boldsymbol{F}$-semimartingale without continuous martingale part. Furthermore, one can easily verify that the $\boldsymbol{F}$-compensator of the jump measure $\mu_{\phi}^{\varepsilon}$ of $\langle \overline{u}^{\varepsilon}, \phi \rangle$ (on $[0,T] \times \mathbb{R}$) is given by
	\begin{equation} \label{compensator nu eps phi}
	\nu_{\phi}^{\varepsilon} (A) = 
	\int_{0}^{T} \int_{0}^{\pi} \int_{\mathbb{R}}  
	\mathbbm{1}_{A} \left( t,\frac{f(u^{\varepsilon}(t,x))}{\sigma(\varepsilon)} \phi(x) z \right)\, \textrm{d}t \, \textrm{d}x \, Q^{\varepsilon}(\textrm{d}z), \quad A \in \mathcal{B}([0,T] \times \mathbb{R}).
	\end{equation}	
	As a consequence, using It\={o}'s formula (see, for example, Theorem I.4.57 in~\cite{JJ}), \eqref{weak formulation} and the fact that $u^{\varepsilon}(0,x) = 0$, we have
	\begin{equation*}
	\begin{split}
	e^{i \xi \langle \overline{u}^{\varepsilon}_t, \phi \rangle} & = 1 + i\xi \int_{0}^{t} e^{i \xi \langle \overline{u}^{\varepsilon}_s, \phi \rangle} \langle u^{\varepsilon}(s, \cdot), \phi'' \rangle \, \textrm{d}s + 
	i \xi \int_{0}^{t} \int_{0}^{\pi} 
	e^{i \xi \langle \overline{u}^{\varepsilon}_{s-}, \phi \rangle} 
	\frac{f(u^{\varepsilon}(s,x))}{\sigma(\varepsilon)}  \phi(x) \,
	L^{\varepsilon}(\textrm{d}s, \textrm{d}x) \\ & \quad \,\, +
	\int_{0}^{t} \int_{\mathbb{R}} 
	e^{i \xi \langle \overline{u}^{\varepsilon}_{s-}, \phi \rangle} 
	( e^{i \xi x} -1 - i \xi x) \, \mu_{\phi}^{\varepsilon}(\textrm{d}s, \textrm{d}x), 
	\end{split}
	\end{equation*}
	and therefore, by \eqref{definition M eps},
	\begin{equation} \label{NR4}
	\begin{split}
	M_t^{\varepsilon} & = 1 + i \xi \int_{0}^{t} \int_{0}^{\pi} 
	e^{i \xi \langle \overline{u}^{\varepsilon}_{s-}, \phi \rangle} 
	\frac{f(u^{\varepsilon}(s,x))}{\sigma(\varepsilon)}  \phi(x) \,
	L^{\varepsilon}(\textrm{d}s, \textrm{d}x) \\ & \quad \,\, + 
	\int_{0}^{t} \int_{\mathbb{R}} 
	e^{i \xi \langle \overline{u}^{\varepsilon}_{s-}, \phi \rangle} 
	(e^{i \xi x} -1 - i \xi x) \, (\mu_{\phi}^{\varepsilon} - \nu_{\phi}^{\varepsilon}) (\textrm{d}s, \textrm{d}x) 
	\end{split}
	\end{equation}
	for all $t \leq T$. The two integral processes on the right-hand side of \eqref{NR4} are square-integrable $\boldsymbol{F}$-martingales (for the second, this is implied by the elementary inequalities $(\cos(x) - 1)^2 \leq x^2$ and $(\sin(x) - x)^2  \leq 4 x^2$ for all $x \in \mathbb{R}$, together with Lemma~\ref{uniform bound second moment}) and hence, this is also the case for $M^{\varepsilon}$. By It\={o}'s isometry, we further obtain
	\begin{equation*}
	\mathbb{E} \left[ \bigg \vert
	\int_{0}^{t}  \int_{\mathbb{R}}  
	e^{i \xi \langle \overline{u}^{\varepsilon}_s, \phi \rangle} 
	\frac{f(u^{\varepsilon}(s,x))}{\sigma(\varepsilon)}  \phi(x) \,
	L^{\varepsilon}(\textrm{d}s, \textrm{d}x)
	\bigg \vert^2
	\right]  = 
	\mathbb{E} \left[
	\int_{0}^{t} \int_{0}^{\pi} 
	f^2(u^{\varepsilon}(s,x)) \phi^2(x) \,
	\textrm{d}s \, \textrm{d}x \right]
	\end{equation*}
	as well as 
	\begin{equation*}
	\begin{split}
	& \mathbb{E} \left[ \bigg \vert \int_{0}^{t}  \int_{\mathbb{R}}  
	e^{i  \xi \langle \overline{u}^{\varepsilon}_s, \phi \rangle} 
	(e^{i  \xi x} -1 - i \xi x) \, (\mu_{\phi}^{\varepsilon} - \nu_{\phi}^{\varepsilon}) (\textrm{d}s, \textrm{d}x)  \bigg \vert^2 \right] \\ & \qquad = 
	\mathbb{E} \left[
	\int_{0}^{t}  \int_{\mathbb{R}}  \vert e^{i  \xi x} -1 - i \xi x \vert^2 \, \nu_{\phi}^{\varepsilon}(\textrm{d}s, \textrm{d}x)
	\right] \\ & \qquad = 
	\mathbb{E} \left[
	\int_{0}^{t}  \int_{\mathbb{R}}  \left( (\cos(\xi x) - 1)^2 + (\sin(\xi x) - \xi x)^2 \right) \, \nu_{\phi}^{\varepsilon}(\textrm{d}s, \textrm{d}x)
	\right]
	\end{split}
	\end{equation*}	
	for all $t \leq T$. We estimate the last expectation, using the elementary inequalities given above as well as the definition of $\nu_{\phi}^{\varepsilon}$, by
	\begin{equation*}
	\begin{split}
	5 \xi^2 \mathbb{E} \left[\int_{0}^{t}  \int_{\mathbb{R}}  x^2 \, \nu_{\phi}^{\varepsilon}(\textrm{d}s, \textrm{d}x) \right] & = 
	5 \xi^2 \mathbb{E} \left[
	\int_{0}^{t} \int_{0}^{\pi}
	\frac{f^2(u^{\varepsilon}(s,x))}{\sigma^2(\varepsilon)} \phi^2(x) z^2 \,
	\textrm{d}s \, \textrm{d}x \, Q^{\varepsilon}(\textrm{d}z) \right]  \\ & = 5  \xi^2
	\mathbb{E} \left[
	\int_{0}^{t} \int_{0}^{\pi} 
	f^2(u^{\varepsilon}(s,x)) \phi^2(x) \,
	\textrm{d}s \, \textrm{d}x \right], \quad t \leq T.
	\end{split}
	\end{equation*}
	Altogether, we obtain~\eqref{bound second moment M eps} from \eqref{NR4}, the Lipschitz continuity of $f$ and Lemma~\ref{Lemma uniform bound}. \qed
\end{Proof}

We now switch to the probability space $(\overline{\Omega}, \overline{\mathcal{F}}, \overline{\boldsymbol{F}}, \overline{\mathbb{P}})$ from the Skorokhod construction in~\eqref{Skorokhod representation 1} and define the process $\overline{M}^k$ in the same way as $M^{\varepsilon}$ in~\eqref{definition M eps}, but with $(u^{\varepsilon}, \overline{u}^{\varepsilon})$ and $Q^{\varepsilon}$ replaced by $(v^k, \overline{v}^k)$ in~\eqref{Skorokhod representation 3} and $Q^{\varepsilon_k}$, respectively.

\begin{Theorem} \label{convergence ov M k to ov M}
	Under~\eqref{AR condition}, we have pointwise on $\overline{\Omega}$ for any $\xi \in \mathbb{R}$ and $\phi \in C_c^{\infty}((0,\pi))$,
	\begin{equation*}
	\overline{M}^k \longrightarrow \overline{M} \quad \textrm{as} \quad k \rightarrow \infty 
	\end{equation*}  
	in the Skorokhod space $D([0,T], \mathbb{C})$, where $\overline{M}$ is the process in~\eqref{def ov M}.
\end{Theorem}	
In order to prove this result, we first rewrite $\overline{M}^k$ in a more convenient form. For any fixed $k \in \mathbb{N}$, $\xi \in \mathbb{R}$ and $\phi \in C_c^{\infty}((0,\pi))$, let
\begin{equation} \label{characteristics of eta k}
\begin{split}
& \overline{\nu}^k (A) = 
\int_{0}^{T} \int_{0}^{\pi} \int_{\mathbb{R}}  
\mathbbm{1}_{A} \left( t, \frac{f(v^k(t,x))}{\sigma(\varepsilon_k)} \phi(x) z \right)  \, \textrm{d}t \, \textrm{d}x \, Q^{\varepsilon_k}(\textrm{d}z),
\\ & \overline{B}_t^{k} = \int_{0}^{t} \langle v^k(s, \cdot), \phi'' \rangle \, \textrm{d}s - 
\int_{0}^{t} \int_{\mathbb{R}} x \mathbbm{1}_{\{ \vert x \vert > 1 \}} \, \overline{\nu}^k(\textrm{d} s, \textrm{d} x), \\
& \overline{A}_t^k = i \xi \overline{B}_t^{k} + \int_{0}^{t} \int_{\mathbb{R}} \left( e^{i \xi x} - 1 - i \xi x \mathbbm{1}_{\{ \vert x \vert \leq 1 \}} \right) \,\overline{\nu}^k(\textrm{d}s, \textrm{d}x)
\end{split}
\end{equation}
for all $A \in \mathcal{B}([0,T] \times \mathbb{R})$ and $t \leq T$. Note that $(\overline{B}^{k}, 0, \overline{\nu}^k)$ are the predictable characteristics of the $\overline{\boldsymbol{F}}$-semimartingale $\langle \overline{v}^k, \phi \rangle$ and that they are functions of the random field $v^k$ and not of $\langle \overline{v}^k, \phi \rangle$ itself (which is another reason why we have adopted a dual view on the solutions to~\eqref{SHE} and~\eqref{SHE 2} as elements of $\Omega^*$). The process $\overline{M}^k$ introduced above can thus be written as
\begin{equation} \label{def ov M k}
\begin{split}
\overline{M}_t^k & = e^{i  \xi \langle \overline{v}^k_t, \phi \rangle}  - 
\int_{0}^{t} e^{i  \xi \langle \overline{v}^k_s, \phi \rangle} \, \overline{A}^k (\textrm{d}s), \quad t \leq T.
\end{split}
\end{equation}

Define the truncation functions
\begin{equation} \label{def truncation function}
\varrho_h \colon \mathbb{R} \longrightarrow \mathbb{R}, \quad x \mapsto  x \mathbbm{1}_{\{ \vert x \vert \leq h \}}, \quad h > 0.
\end{equation}
The key idea of the proof of Theorem~\ref{convergence ov M k to ov M} is to see that for fixed $\overline{\omega} \in \overline{\Omega}$, $t \leq T$ and $\phi \in C_c^{\infty}((0,\pi))$, the function $\overline{A}_t^k$ in~\eqref{characteristics of eta k} (resp., $\overline{A}_t$ in~\eqref{characteristics of eta}) is the Lévy exponent of the infinitely divisible distribution $\eta_k$ (resp., $\eta$) with characteristics $(\overline{B}_t^k, 0, \overline{\nu}^k([0,t] \times \textrm{d}x))$ (resp., $(\overline{B}_t, \overline{C}_t, 0)$) with respect to $\varrho_1$. Then we can make use of the following result, which is the only place in this work where \eqref{AR condition} will actually be needed.

%
%
%

\begin{Theorem} \label{weak convergnce id distributions}
	If \eqref{AR condition} holds, then for any $\phi \in C_c^{\infty}((0,\pi))$, $t \in [0,T]$ and $\overline{\omega} \in \overline{\Omega}$, we have
	\begin{equation*}
	\eta_k \stackrel{w}{\longrightarrow} \eta \quad \textrm{as} \quad k \rightarrow \infty.
	\end{equation*}
\end{Theorem}

\begin{Proof}
As in the proof of Theorem~2.2 in~\cite{CR} (see also Theorem~2.1 in \cite{AR}), it suffices to show that
\begin{equation} \label{convergence of deterministic characteristics}
\begin{split}
& (i) \quad \int_{\vert x \vert \leq h}  x^2 \, \overline{\nu}^k([0,t] \times \textrm{d}x) \longrightarrow \overline{C}_t, \\ 
& (ii) \quad \overline{B}_t^k \longrightarrow \overline{B}_t,   \\ 
& (iii) \quad \overline{\nu}^k \left([0,t] \times \{ \vert x \vert > 1 \}\right) \longrightarrow 0
\end{split}
\end{equation}
as $k \rightarrow \infty$ for all $h > 0$. Starting with $(i)$, we have
\begin{equation*}
\begin{split}
& \int_{\vert x \vert \leq h} x^2 \, \overline{\nu}^k([0,t] \times \textrm{d}x) \\ & \qquad = 
\int_{0}^{t} \int_{0}^{\pi} 
f^2(v^k(s,x)) \phi^2(x) \frac{1}{\sigma^2(\varepsilon_k)} \int_{\mathbb{R}} z^2 
\mathbbm{1}_{ \left \{ \vert z \vert \leq \left(h/\vert f(v^k(s,x)) \phi(x) \vert\right) \sigma(\varepsilon_k) \right \} } \,
Q^{\varepsilon_k}(\textrm{d}z) \, \textrm{d}s \, \textrm{d}x.
\end{split}
\end{equation*}
We can ignore all points in the domain of integration where $\vert f(v^k(s,x)) \phi(x) \vert= 0$. So if we let
\begin{equation} \label{def I and Sigma}
\begin{split}
I^k_h(s,x) & = \frac{1}{\sigma^2(\varepsilon_k)}  \int_{\mathbb{R}} z^2 
\mathbbm{1}_{ \left \{ \vert z \vert \leq \left(h/\vert f(v^k(s,x)) \phi(x) \vert\right)  \sigma(\varepsilon_k) \right \} } \,
Q^{\varepsilon_k}(\textrm{d}z) \quad \textrm{and} \\
\Sigma^k_h(s,x) & = 1 - I^k_h(s,x) = 
\frac{1}{\sigma^2(\varepsilon_k)} 
\int_{\mathbb{R}} z^2  
\mathbbm{1}_{ \left \{ \vert z \vert > \left(h/\vert f(v^k(s,x)) \phi(x) \vert\right) \sigma(\varepsilon_k) \right \}} \, Q^{\varepsilon_k}(\textrm{d}z)
\end{split}
\end{equation} 
for any $(s,x) \in [0,T] \times [0,\pi]$, $h > 0$ and $k \in \mathbb{N}$, then, using the triangle inequality and the fact that $0 < I^k_h(s,x) \leq 1$, we obtain 
\begin{equation} \label{NR1}
\begin{split}
& \Bigg \vert \int_{\vert x \vert \leq h} x^2 \, \overline{\nu}^k([0,t] \times \textrm{d}x) - \overline{C}_t \Bigg \vert \\ & \qquad \leq
\Bigg \vert 
\int_{0}^{t} \int_{0}^{\pi}
\left(f^2(v^k(s,x)) - f^2(v(s,x))\right) \phi^2(x) 
I^k_h(s,x)
\, \textrm{d}s \, \textrm{d}x \Bigg \vert \\ & \qquad  \quad \,\,  + \Bigg \vert
\int_{0}^{t} \int_{0}^{\pi}
f^2(v(s,x)) \phi^2(x)  \left( I^k_h(s,x) - 1 \right)\,
\textrm{d}s \, \textrm{d}x \Bigg \vert \\ & \qquad  \leq 
\int_{0}^{t} \int_{0}^{\pi}
\Big \vert f^2(v^k(s,x)) - f^2(v(s,x)) \Big \vert  \phi^2(x) \,\textrm{d}s \, \textrm{d}x +
\int_{0}^{t} \int_{0}^{\pi}
f^2(v(s,x)) \phi^2(x)  \Sigma^k_h(s,x) \, \textrm{d}s \, \textrm{d}x.
\end{split}
\end{equation}
By the Lipschitz continuity of $f$ and Hölder's inequality, we have for the first integral on the right-hand side of \eqref{NR1},
\begin{equation} \label{NR10}
\begin{split}
& \int_{0}^{t} \int_{0}^{\pi}
\Big \vert f^2(v^k(s,x)) - f^2(v(s,x)) \Big \vert  \phi^2(x) \, \textrm{d}s \, \textrm{d}x \\ & \qquad \leq 
C \int_{0}^{t} \int_{0}^{\pi} 
\Big \vert f(v^k(s,x)) - f(v(s,x)) \Big \vert 
\Big \vert f(v^k(s,x)) + f(v(s,x)) \Big \vert \,\textrm{d}s \, \textrm{d}x \\ & \qquad \leq
C \left(\int_{0}^{t} \int_{0}^{\pi} 
\left( v^k(s,x) - v(s,x) \right)^2 \,\textrm{d}s \, \textrm{d}x\right)^{1/2} \left(\int_{0}^{t} \int_{0}^{\pi} 
\left( f(v^k(s,x)) + f(v(s,x)) \right)^2 \,\textrm{d}s \, \textrm{d}x\right)^{1/2}.
\end{split}
\end{equation}
By~\eqref{Skorokhod representation 3}, $v^k \longrightarrow v$ in $L^2([0,T] \times [0,\pi])$ pointwise on $\overline{\Omega}$. Hence, the sequence ${(v^k)}_{k \in \mathbb{N}}$ is bounded in $L^2([0,T] \times [0,\pi])$ and we have
\begin{equation} \label{NR11}
\begin{split}
& \int_{0}^{t} \int_{0}^{\pi} 
\left( f(v^k(s,x)) + f(v(s,x)) \right)^2 \, \textrm{d}s \, \textrm{d}x  \\ & \qquad \leq
C \left(1 + \sup_{k \in \mathbb{N}} \int_{0}^{t} \int_{0}^{\pi} v^k(s,x)^2 \, \textrm{d}s \, \textrm{d}x  + 
\int_{0}^{t} \int_{0}^{\pi} v(s,x)^2 \, \textrm{d}s \, \textrm{d}x \right)  < \infty,
\end{split}
\end{equation}
which implies 
$\int_{0}^{t} \int_{0}^{\pi}
\vert f^2(v^k(s,x)) - f^2(v(s,x)) \vert  \phi^2(x) \, \textrm{d}s \, \textrm{d}x \longrightarrow 0$ as $k \rightarrow \infty$.

The second integral in \eqref{NR1} is more difficult to handle. We decompose it into $I_1^{k,n} + I_2^{k,n,M} + I_3^{k,n,M}$, where
\begin{equation*}
\begin{split}
I_1^{k,n} & =
\int_{0}^{t} \int_{0}^{\pi}
f^2(v(s,x)) 
\Sigma^k_h(s,x) \mathbbm{1}_{\{ \vert v^k(s,x) \vert \leq n \}} \, \textrm{d}s \, \textrm{d}x, \\ 
I_2^{k,n,M} & =
\int_{0}^{t} \int_{0}^{\pi}
f^2(v(s,x)) 
\Sigma^k_h(s,x) \mathbbm{1}_{\{ \vert v^k(s,x) \vert > n \}} 
\mathbbm{1}_{\{ \vert f(v(s,x)) \vert \leq M \}} \, \textrm{d}s \, \textrm{d}x, \\ 
I_3^{k,n,M} & =
\int_{0}^{t} \int_{0}^{\pi}
f^2(v(s,x))  
\Sigma^k_h(s,x)  \mathbbm{1}_{\{ \vert v^k(s,x) \vert > n \}}  
\mathbbm{1}_{\{ \vert f(v(s,x)) \vert > M \}} \, \textrm{d}s \, \textrm{d}x
\end{split}
\end{equation*}
for all $k, n, M \in \mathbb{N}$. 

Again we will study each of these three integrals separately. On the set $\{ \vert v^k(s,x) \vert \leq n \}$, we have $\vert f(v^k(s,x)) \phi(x) \vert \leq (Kn + \vert f(0) \vert ) {\Vert \phi \Vert}_{\infty}$ and therefore
\begin{equation*}
\mathbbm{1}_{ \left \{ \vert z \vert > \left(h/\vert f(v^k(s,x)) \phi(x) \vert\right) \sigma(\varepsilon_k) \right \}}  \leq 
\mathbbm{1}_{ \left \{ \vert z \vert > \left(h/(Kn + \vert f(0) \vert ) {\Vert \phi \Vert}_{\infty}\right) \sigma(\varepsilon_k) \right \}}.
\end{equation*}
Thus,
\begin{equation*}
I_1^{k,n} \leq
\int_{0}^{t} \int_{0}^{\pi}
f^2(v(s,x))  
\frac{1}{\sigma^2(\varepsilon_k)} 
\int_{\mathbb{R}} z^2 \mathbbm{1}_{ \left \{ \vert z \vert > \left(h/(Kn + \vert f(0) \vert ) {\Vert \phi \Vert}_{\infty}\right) \sigma(\varepsilon_k) \right \}} \, Q^{\varepsilon_k}(\textrm{d}z) \, \textrm{d}s \, \textrm{d}x.
\end{equation*}
Because the term $h/ (Kn + \vert f(0) \vert ) {\Vert \phi \Vert}_{\infty}$ does not depend on $k$, we can use condition~\eqref{AR condition}, whence
\begin{equation*}
\frac{1}{\sigma(\varepsilon_k)}  
\int_{\mathbb{R}} z^2  \mathbbm{1}_{ \left \{ \vert z \vert > \left(h/(Kn + \vert f(0) \vert ) {\Vert \phi \Vert}_{\infty}\right) \sigma(\varepsilon_k) \right \}} \, Q^{\varepsilon_k}(\textrm{d}z)
\longrightarrow 0 \quad \textrm{as} \quad k \rightarrow \infty
\end{equation*}
for all $n \in \mathbb{N}$ and $h > 0$. Since $v \in L^2([0,T] \times [0,\pi])$, we obtain by dominated convergence that $I_1^{k,n} \longrightarrow 0$ as $k \rightarrow \infty$ for all $n \in \mathbb{N}$. 

Next, we have by Chebyshev's inequality,
\begin{equation*}
I_2^{k,n,M} \leq
M^2 \int_{0}^{t} \int_{0}^{\pi} \mathbbm{1}_{\{ \vert v^k(s,x) \vert > n \}} \, \textrm{d}s \, \textrm{d}x \leq 
\frac{M^2}{n^2} 
\sup_{k \in \mathbb{N}} 
\int_{0}^{t} \int_{0}^{\pi} v^k(s,x)^2 \, \textrm{d}s \, \textrm{d}x,
\end{equation*}
which tends to 0 as $n \rightarrow \infty$, uniformly in $k$.

Finally, we have by dominated convergence,
\begin{equation*}
I_3^{k,n,M} \leq
\int_{0}^{t} \int_{0}^{\pi}
f^2(v(s,x))  
\mathbbm{1}_{\{ \vert f(v(s,x)) \vert > M \}} \, \textrm{d}s \, \textrm{d}x \longrightarrow 0 \quad \textrm{as} \quad M \rightarrow \infty,
\end{equation*}
uniformly in $n$ and $k$. Altogether, we have just shown that the left-hand side of \eqref{NR1} converges to 0 as $k \rightarrow \infty$, which is condition $(i)$ in~\eqref{convergence of deterministic characteristics}. 

The two other conditions will follow from our last calculations. Indeed, we have 
\begin{equation*}
\big \vert 
\overline{B}_t^k - \overline{B}_t
\big \vert \leq
\int_{0}^{t} \Big \vert \langle v^k(s, \cdot), \phi'' \rangle - \langle v(s,\cdot), \phi'' \rangle \Big \vert \, \textrm{d}s + 
\Bigg \vert \int_{0}^{t} \int_{\mathbb{R}} x \mathbbm{1}_{\{ \vert x \vert > 1 \}} \, \overline{\nu}^k(\textrm{d} s, \textrm{d} x) \Bigg \vert, 
\end{equation*}
where the first term vanishes because
\begin{equation} \label{NR15}
\begin{split}
\int_{0}^{t} \Big \vert \langle v^k_s, \phi'' \rangle - \langle v_s, \phi'' \rangle \Big \vert \, \textrm{d}s & \leq
\int_{0}^{t} \int_{0}^{\pi} \Big \vert v^k(s,x) - v(s,x) \Big \vert  \vert \phi''(x)  \vert \, \textrm{d}s \, \textrm{d}x \\ & \leq
C \left(\int_{0}^{t} \int_{0}^{\pi} \left( v^k(s,x) - v(s,x)\right)^2 \,\textrm{d}s \, \textrm{d}x \right)^{1/2}
\longrightarrow 0
\end{split}
\end{equation}
as $k \rightarrow \infty$. Furthermore,
\begin{equation} \label{NR16}
\begin{split}
 \int_{0}^{t} \int_{\mathbb{R}} x \mathbbm{1}_{\{ \vert x \vert > 1 \}} \, \overline{\nu}^k(\textrm{d} s, \textrm{d} x)  &\leq 
\int_{0}^{t} \int_{\mathbb{R}} x^2 \mathbbm{1}_{\{ \vert x \vert > 1 \}} \, \overline{\nu}^k(\textrm{d} s, \textrm{d} x) \\ &  = 
\int_{0}^{t} \int_{0}^{\pi} 
f^2(v^k(s,x)) \phi^2(x) \Sigma^k_1(s,x) \, \textrm{d}s \, \textrm{d}x \\ &  \leq
\int_{0}^{t} \int_{0}^{\pi}  \big \vert f^2(v^k(s,x)) - f^2(v(s,x)) \big \vert \phi^2(x) \, \textrm{d}s \, \textrm{d}x \\ &  \quad \,\,  +
\int_{0}^{t} \int_{0}^{\pi} 
f^2(v(s,x)) \phi^2(x) \Sigma^k_1(s,x) \, \textrm{d}s \, \textrm{d}x.
\end{split}
\end{equation}
These integrals are exactly the same as in the last line of (\ref{NR1}), so we obtain 
$\big \vert 
\overline{B}_t^k - \overline{B}_t
\big \vert \longrightarrow 0$, which is condition $(ii)$. From this, condition $(iii)$ immediately follows since
\begin{equation} \label{NR14}
\overline{\nu}^k \left([0,t] \times \{ \vert x \vert > 1 \}\right) = 
\int_{0}^{t} \int_{\mathbb{R}} \mathbbm{1}_{\{ \vert x \vert > 1 \}} \, \overline{\nu}^k(\textrm{d} s, \textrm{d} x)  \leq
\int_{0}^{t} \int_{\mathbb{R}} x \mathbbm{1}_{\{ \vert x \vert > 1 \}} \, \overline{\nu}^k(\textrm{d} s, \textrm{d} x).
\end{equation}
\qed
\end{Proof}

The following technical lemma is a direct consequence of Theorem~\ref{weak convergnce id distributions} and will be crucial for proving Theorem~\ref{convergence ov M k to ov M} afterwards.

If $t \mapsto A_t$ is a function of locally finite variation, we denote by  $\textrm{Var}(A)_t$ the total variation of the function $A$ on the interval $[0,t]$. If $A$ is complex-valued, we have $\textrm{Var}(A) = \textrm{Var} (\mathrm{Re} \, A) + \textrm{Var} (\mathrm{Im} \, A)$.

\begin{Lemma} \label{convergence of the variation}
	If~\eqref{AR condition} holds, then we have pointwise on $\overline{\Omega}$, 
	\begin{equation*}
	\overline{A}^k_t \longrightarrow \overline{A}_t \quad \textrm{and} \quad \emph{\textrm{Var}}(\overline{A}^k - \overline{A})_t \longrightarrow 0 \quad \textrm{as} \quad k \rightarrow \infty
	\end{equation*}	
	for any $\xi \in \mathbb{R}$, $\phi \in C_c^{\infty}((0,\pi))$ and $t \leq T$, where the processes $\overline{A}^k$ and $\overline{A}$ are defined in~\eqref{characteristics of eta k} and~\eqref{characteristics of eta}, respectively.
\end{Lemma}

\begin{Proof}

For fixed $\phi \in C_c^{\infty}((0,\pi))$, $t \in [0,T]$ and $\overline{\omega} \in \overline{\Omega}$, the infinitely divisible distributions $\eta_k$ and $\eta$, defined before Theorem~\ref{weak convergnce id distributions}, have Lévy exponents $\overline{A}_t^k$ and $\overline{A}_t$, respectively. By that theorem, $\eta_k \stackrel{w}{\longrightarrow} \eta$ as $k \rightarrow \infty$. This immediately implies the first claim of the proposition (see, for example, Equation VII.2.6 in \cite{JJ}).

For the second claim, we will need the truncation function
\begin{equation*}
\vartheta(x) = \left\{ 
\begin{array}{cc}
-1, & x < -1, \\
x, & \vert x \vert \leq 1, \\
1, & x > 1.
\end{array} \right.
\end{equation*}
The main difference between the function $\varrho_1(x) = x \mathbbm{1}_{\{ \vert x\vert \leq 1\}}$, used so far, and $\vartheta$ is that the latter is continuous. Since this property will be needed for technical reasons, we replace $\varrho_1$ by $\vartheta$ in the expression of $\overline{A}_t^k$ in~\eqref{characteristics of eta k} and thus obtain
\begin{equation} \label{NR25}
\begin{split}
\overline{A}_t^k & =
i \xi  \left(\int_{0}^{t} \langle v^k(s, \cdot), \phi'' \rangle \, \textrm{d}s - 
\int_{0}^{t} \int_{\mathbb{R}}  \left(x - \vartheta(x)\right) \,\overline{\nu}^k(\textrm{d} s, \textrm{d} x) \right) \\ & \quad \,\, + \int_{0}^{t}  \int_{\mathbb{R}} \left( e^{i \xi x} - 1 - i \xi \vartheta(x) \right) \,\overline{\nu}^k(\textrm{d} s, \textrm{d} x).
\end{split}
\end{equation}
With~\eqref{NR25} and~\eqref{characteristics of eta}, we then calculate
\begin{equation} \label{NR3}
\begin{split}
\textrm{Re} (\overline{A}^k_t - \overline{A}_t) & =
\frac{1}{2} \xi^2 \int_{0}^{t} \int_{0}^{\pi} f^2(v(s,x)) \phi^2(x) \, \textrm{d}s \, \textrm{d}x +
\int_{\mathbb{R}} \left( \cos (\xi x) -1 \right)\, \overline{\nu}^k([0,t] \times \textrm{d}x), \\ 
\textrm{Im} (\overline{A}^k_t - \overline{A}_t) & =
\xi \left( \int_{0}^{t} \langle v^k(s, \cdot), \phi'' \rangle \, \textrm{d}s - \int_{0}^{t} \langle v(s, \cdot), \phi'' \rangle \, \textrm{d}s - 
\int_{\mathbb{R}} \left(x - \vartheta(x)\right)\, \overline{\nu}^k([0,t] \times \textrm{d}x) \right) \\ & \quad \,\, + 
\int_{\mathbb{R}} \left( \sin (\xi x) - \xi \vartheta(x) \right)\, \overline{\nu}^k([0,t] \times \textrm{d}x).
\end{split}
\end{equation}
Consequently,
\begin{equation} \label{NR2}
\begin{split}
\textrm{Var}(\textrm{Re} (\overline{A}^k - \overline{A}))_t &  \leq 
\frac{1}{2} \xi^2 \int_{0}^{t} 
\int_{0}^{\pi} \Bigg \vert f^2(v(s,x)) \phi^2(x)  - 
\int_{\mathbb{R}} \vartheta^2  \left( \frac{f(v^k(s,x))}{\sigma(\varepsilon_k)} \phi(x) z \right) Q^{\varepsilon_k}(\textrm{d}z) \Bigg \vert \, \textrm{d}x \, \textrm{d}s \\ & \quad \,\, + 
\int_{\mathbb{R}} \Big \vert \cos \left( \xi x \right) - 1  +
\frac{1}{2} \xi^2 \vartheta^2 (x) \Big \vert \, \overline{\nu}^k([0,t] \times \textrm{d}x).
\end{split}
\end{equation}
We will show that the two integrals above converge to 0 as $k \rightarrow \infty$. 

The function $\vert \cos(\xi x) - 1 + \frac{1}{2} \xi^2 \vartheta^2(x) \vert$ is bounded, continuous (because $\vartheta$ is) and $o(x^2)$ as $x \rightarrow 0$. Hence, by condition $[\delta_{1,3}]$ of Theorem VII.2.9 in \cite{JJ}, we can infer
\begin{equation*}
\int_{\mathbb{R}} \Big \vert \cos(\xi x) - 1 + \frac{1}{2} \xi^2 \vartheta^2(x) \Big \vert \, \overline{\nu}^k([0,t] \times \textrm{d} x) \longrightarrow 0 \quad \textrm{as} \quad k \rightarrow \infty
\end{equation*}
from Theorem \ref{weak convergnce id distributions}.

Consider now the first integral in (\ref{NR2}), and notice that 
\begin{equation*}
\vartheta^2(x) = x^2 \mathbbm{1}_{\{ \vert x \vert \leq 1 \}} + \mathbbm{1}_{\{ \vert x \vert > 1 \}} \quad \textrm{for all} \quad x \in \mathbb{R}.
\end{equation*}
Using the triangle inequality, we can therefore estimate this integral by
\begin{equation*}
\begin{split}
& \int_{0}^{t} \int_{0}^{\pi} \Bigg \vert
f^2(v(s,x)) \phi^2(x) - 
\int_{\mathbb{R}} \left( f(v^k(s,x)) \phi(x) z /\sigma(\varepsilon_k)  \right)^2 
\mathbbm{1}_{\left \{ \vert f(v^k(s,x)) \phi(x) z \vert /\sigma(\varepsilon_k) \leq 1 \right \}} \,Q^{\varepsilon_k}(\textrm{d}z) \Bigg \vert \, \textrm{d}x \, \textrm{d}s \\ & \qquad +
\int_{\mathbb{R}} \mathbbm{1}_{\{ \vert x \vert > 1 \}} \, \overline{\nu}^k([0,t] \times \textrm{d} x).
\end{split}
\end{equation*}
The second integral above converges to 0 as shown in \eqref{NR14}, while the same holds for the first integral by \eqref{NR1} (set $h = 1$). Together with (\ref{NR2}), we conclude that $\textrm{Var}(\textrm{Re} (\overline{A}^k - \overline{A}))_t \longrightarrow 0$ as $k \rightarrow \infty$.

It remains to show that $\textrm{Var}(\textrm{Im} (\overline{A}^k - \overline{A}))_t \longrightarrow 0$ as $k \rightarrow \infty$, which will be done in a similar manner as before. From (\ref{NR3}), we have
\begin{equation*}
\begin{split}
\textrm{Var}(\textrm{Im} (\overline{A}^k - \overline{A}))_t  & \leq
\vert \xi \vert \int_{0}^{t} \Big \vert 
\langle v^k(s, \cdot), \phi'' \rangle - \langle v(s, \cdot), \phi'' \rangle \Big \vert \, \textrm{d}s + \vert \xi \vert 
\int_{\mathbb{R}} \big \vert x - \vartheta(x) \big \vert \, \overline{\nu}^k([0,t] \times \textrm{d}x) \\ & \quad \,\, +
\int_{\mathbb{R}} \Big \vert \sin (\xi x) - \xi \vartheta(x) \Big \vert \, \overline{\nu}^k([0,t] \times \textrm{d}x).
\end{split}
\end{equation*}
The first integral on the right-hand side above converges to 0 by \eqref{NR15}. Furthermore, since
\begin{equation*}
\int_{\mathbb{R}} \big \vert x - \vartheta(x) \big \vert \, \overline{\nu}^k([0,t] \times \textrm{d}x) \leq 
2 \int_{\mathbb{R}} {\vert x \vert}^2 \mathbbm{1}_{\{ x > 1 \}} \, \overline{\nu}^k([0,t] \times \textrm{d}x),
\end{equation*}
also the second integral vanishes by \eqref{NR16}. Finally, the function $\big \vert \sin (\xi x) - \xi \vartheta(x) \big \vert$ is bounded, continuous and  $o(x^2)$ as $x \rightarrow 0$. Hence, we can again apply Theorem VII.2.9 in \cite{JJ} in order to obtain
\begin{equation*}
\int_{\mathbb{R}} \Big \vert \sin (\xi x) - \xi \vartheta(x) \Big \vert \, \overline{\nu}^k([0,t] \times \textrm{d}x) \longrightarrow 0 \quad \textrm{as} \quad k \rightarrow \infty.
\end{equation*}
We conclude that $\textrm{Var}(\textrm{Im} (\overline{A}^k - \overline{A}))_t \longrightarrow 0$, and altogether $\textrm{Var}(\overline{A}^k - \overline{A})_t \longrightarrow 0$ as $k \rightarrow \infty$. \qed 
\end{Proof}

\begin{Proof}[of Theorem~\ref{convergence ov M k to ov M}]
According to Proposition VI.1.23 in \cite{JJ}, because the function
$t \mapsto \int_{0}^{t} \exp \left(i  \xi \langle \overline{v}_s, \phi \rangle \right) \,\overline{A} (\textrm{d}s)$ is continuous, $\overline{M}^k$ converges to $\overline{M}$ in the Skorokhod topology for fixed $\overline{\omega} \in \overline{\Omega}$ if 
\begin{equation*}
e^{i  \xi \langle \overline{v}^k, \phi \rangle } \longrightarrow e^{i  \xi \langle \overline{v}, \phi \rangle }
\quad \textrm{and} \quad
\int_{0}^{\cdot} e^{i  \xi \langle \overline{v}^k_s, \phi \rangle } \, \overline{A}^k (\textrm{d}s) \longrightarrow 
\int_{0}^{\cdot} e^{i  \xi \langle \overline{v}_s, \phi \rangle } \, \overline{A} (\textrm{d}s)
\end{equation*}
in $D([0,T], \mathbb{C})$ as $k \rightarrow \infty$.

Using the definition of the Skorokhod topology, we can easily infer from the convergence of ${(\overline{v}^k)}_{k \in \mathbb{N}}$ to $\overline{v}$ in $D([0,T], H_{-r}([0,\pi]))$ given in~\eqref{Skorokhod representation 3} that
\begin{equation*}
\begin{split}
& \langle \overline{v}^k, \phi \rangle \longrightarrow \langle \overline{v}, \phi \rangle \quad \textrm{in} \quad D([0,T], \mathbb{R}) \quad \textrm{and} \\
& e^{i \xi \langle \overline{v}^k, \phi \rangle} \longrightarrow e^{i \xi \langle \overline{v}, \phi \rangle} \quad \textrm{in} \quad D([0,T], \mathbb{C})
\end{split}
\end{equation*}
as $k \rightarrow \infty$ for all $\phi \in H_r([0, \pi])$.

Next, we have
\begin{equation*}
\begin{split}
& \sup_{t \leq T} \Bigg \vert \int_{0}^{t} e^{i  \xi \langle \overline{v}^k_s, \phi \rangle} \, \overline{A}^k (\textrm{d}s) - \int_{0}^{t} e^{i \xi \langle \overline{v}_s, \phi \rangle} \, \overline{A} (\textrm{d}s) \Bigg \vert \\ & \qquad \leq
\sup_{t \leq T} \Bigg \vert 
\int_{0}^{t} e^{i  \xi \langle \overline{v}^k_s, \phi \rangle} \, (\overline{A}^k - \overline{A} ) (\textrm{d}s)  \Bigg \vert +
\sup_{t \leq T} \Bigg \vert 
\int_{0}^{t} \left( e^{i  \xi \langle \overline{v}^k_s, \phi \rangle} - e^{i \xi \langle \overline{v}_s, \phi \rangle} \right)\, \overline{A} (\textrm{d}s) \Bigg \vert \\ & \qquad  \leq
{\textrm{Var} (\overline{A}^k - \overline{A} )}_T  +
\int_{0}^{T} \big \vert e^{i  \xi \langle \overline{v}^k_s, \phi \rangle} - e^{i  \xi \langle \overline{v}_s, \phi \rangle} \big \vert \, \textrm{Var}(\overline{A}) (\textrm{d}s).
\end{split}
\end{equation*}
Lemma \ref{convergence of the variation} then immediately gives us ${\textrm{Var} (\overline{A}^k - \overline{A} )}_T \rightarrow 0$ as $k \rightarrow \infty$. In addition, the Skorokhod convergence of $e^{i \xi \langle \overline{v}^k, \phi \rangle}$ towards $e^{i \xi \langle \overline{v}, \phi \rangle}$ implies 
$e^{i \xi \langle \overline{v}^k_t, \phi \rangle} \longrightarrow e^{i \xi \langle \overline{v}_t, \phi \rangle}$ for all continuity points of $e^{i \xi \langle \overline{v}, \phi \rangle}$; see, for example, VI.2.3 of~\cite{JJ}. Since a càdlàg function has at most countably many discontinuities, we have $e^{i \xi \langle \overline{v}^k_t, \phi \rangle} \longrightarrow e^{i \xi \langle \overline{v}_t, \phi \rangle}$ for almost all $t \in [0,T]$. So dominated convergence implies that also the last term of the previous display converges to 0 as $k \rightarrow \infty$. \qed 
\end{Proof}

We have now gathered all the intermediate results needed for the following theorem.

\begin{Theorem} \label{convergence of the MPs}	
	If~\eqref{AR condition} holds, then $(v, \overline{v})$ in~\eqref{Skorokhod representation 3} satisfies the following martingale problem. For all $\xi \in \mathbb{R}$ and $\phi \in C_c^{\infty}((0,\pi))$, the process $(\overline{M}_t)_{t \leq T}$ defined in~\eqref{def ov M} is a martingale with respect to the filtration $\overline{\boldsymbol{F}}$ in~\eqref{def filtration Om Stern}.
		
	Furthermore, $v$ has an $\overline{\boldsymbol{F}}$-predictable modification and
	\begin{equation} \label{ess sup 2}
	\esssup_{(t,x) \in [0,T] \times [0,\pi]}  \mathbb{E}\left[ \vert v(t,x) \vert^2 \right] < \infty.
	\end{equation}
	Finally, for almost all $t \in [0,T]$, $\overline{v}_t = \langle v(t, \cdot), \cdot \rangle$ as well as $\overline{v}_0=0$ holds with probability one. 
\end{Theorem}

\begin{Proof}
By Theorem~\ref{solution of SHE with Lévy noise is solution of MP}, for any $\xi \in \mathbb{R}$, $\phi \in C_c^{\infty}((0, \pi))$ and $k \in \mathbb{N}$, the process $M^{\varepsilon_k}$ defined in~\eqref{definition M eps} is a square-integrable $\boldsymbol{F}$-martingale. Moreover, as $\overline{v}^k$ and $v^k$ in~\eqref{Skorokhod representation 3} are adapted to the filtration $\overline{\boldsymbol{F}}$, the same holds for $\overline{M}^k$ from~\eqref{def ov M k} as well as $\overline{v}$, $v$ and $\overline{M}$ by a limit argument. Since $\overline{M}^k$ has the same distribution as $M^{\varepsilon_{k}}$ by~\eqref{Skorokhod representation 3}, standard arguments now show that $\overline{M}^k$ is an $\overline{\boldsymbol{F}}$-martingale for all $\xi \in \mathbb{R}$, $\phi \in C_c^{\infty}((0, \pi))$ and $k \in \mathbb{N}$. This is the martingale problem satisfied by the pair $(v^{k}, \overline{v}^k)$. 

Using Theorem~\ref{convergence ov M k to ov M}, we have
\begin{equation} \label{a.s. convergence in Skorokhod space}
\overline{M}^k(\overline{\omega}) \longrightarrow \overline{M}(\overline{\omega}) \quad \textrm{in} \quad D([0,T], \mathbb{C})  
\end{equation}  
as $k \rightarrow \infty$ for all $\overline{\omega} \in \overline{\Omega}$. This implies $\overline{M}^k_t(\overline{\omega}) \longrightarrow \overline{M}_t(\overline{\omega})$ almost everywhere on $[0,T]$ for all $\overline{\omega} \in \overline{\Omega}$. Furthermore, 
\begin{equation*}
\mathbb{E}\left[ \big \vert \overline{M}^k_t \big \vert^2 \right] = 
\mathbb{E}\left[ \big \vert M^{\varepsilon_{k}}_t \big \vert^2 \right] < \infty
\end{equation*}
uniformly in $k \in \mathbb{N}$ and $t \leq T$ by Theorem~\ref{solution of SHE with Lévy noise is solution of MP}. Hence, again by standard arguments, we can deduce that $\overline{M}$ is an $\overline{\boldsymbol{F}}$-martingale as well for any $\xi \in \mathbb{R}$ and $\phi \in C_c^{\infty}((0, \pi))$. 

Now we show the second part of the theorem. The convergence in~\eqref{Skorokhod representation 3} implies convergence in measure (with respect to the Lebesgue measure on $[0,T] \times [0,\pi]$) of $v^k(\overline{\omega})$ towards $v(\overline{\omega})$ for all $\overline{\omega} \in \overline{\Omega}$. Hence, we have by dominated convergence,
\begin{equation*}
\overline{\mathbb{P}} \otimes {\textrm{Leb}}_{[0,T] \times [0,\pi]} \left( \vert v^k - v \vert \geq \varepsilon \right) = 
\mathbb{E} \left[ \int_{0}^{T} \int_{0}^{\pi} \mathbbm{1}_{\left \{ \vert v^k(\overline{\omega}, t,x) - v(\overline{\omega}, t,x) \vert \geq \varepsilon \right \}} \, \textrm{d}t \, \textrm{d}x \right] \longrightarrow 0
\end{equation*}
as $k \rightarrow \infty$, and thus, $v^k$ converges to $v$ in $\overline{\mathbb{P}} \otimes {\textrm{Leb}}_{[0,T] \times [0,\pi]}$-measure. Therefore, there exists a subsequence $(k_l)_{l \in \mathbb{N}}$ such that 
\begin{equation} \label{pointwise convergence of v k}
v^{k_l} \longrightarrow v \quad \overline{\mathbb{P}} \otimes {\textrm{Leb}}_{[0,T] \times [0,\pi]} \textrm{-almost everywhere} 
\quad  \textrm{as} \quad l \rightarrow \infty,
\end{equation}
and we will assume without loss of generality that~\eqref{pointwise convergence of v k} holds for the whole sequence. This in turn implies
$v^{k} \longrightarrow v$ $\overline{\mathbb{P}}$-almost surely as $k \rightarrow \infty$ for almost all $(t,x) \in [0,T] \times [0,\pi]$. Using Fatou's lemma, we obtain 
\begin{equation} \label{Fatou lemma}
\mathbb{E} \left[ \vert v(t,x) \vert^2 \right] \leq \liminf_{k \rightarrow \infty} \mathbb{E} \left[ \vert v^{k}(t,x) \vert^2 \right] \quad 
{\textrm{Leb}}_{[0,T] \times [0,\pi]} \textrm{-almost everywhere}.
\end{equation}
Furthermore, 
\begin{equation} \label{a.e. equality in distribution v k and u eps k}
v^k(t,x) \stackrel{d}{=} u^{\varepsilon_k}(t,x) \quad {\textrm{Leb}}_{[0,T] \times [0,\pi]} \textrm{-almost everywhere},
\end{equation}
for all $k \in \mathbb{N}$, so~\eqref{ess sup 2} follows from Lemma~\ref{uniform bound second moment}. (In order to show~\eqref{a.e. equality in distribution v k and u eps k}, consider for $\alpha > 0$ the mollified random fields $J_{\alpha} v^k$ and $J_{\alpha} u^{\varepsilon_k}$ on $[0,T] \times [0,\pi]$, defined exactly as in~(1.8) of Chapter~10 in~\cite{Friedman}. Then~\eqref{Skorokhod representation 3} implies 
\begin{equation} \label{mollifier 1}
(J_{\alpha} v^k)(t,x) \stackrel{d}{=} (J_{\alpha} u^{\varepsilon_k})(t,x) 
\end{equation}
for all $(t,x) \in [0,T] \times [0,\pi]$, $\alpha > 0$ and $k \in \mathbb{N}$. In addition, using Lemma~3 of Chapter~10 in~\cite{Friedman}, we have
\begin{equation*}
J_{\alpha} v^k(\overline{\omega}) \longrightarrow v^k(\overline{\omega}) \quad \textrm{and} \quad J_{\alpha} u^{\varepsilon_k}(\omega) \longrightarrow u^{\varepsilon_k}(\omega) \quad \textrm{in} \quad L^2([0,T] \times [0,\pi]) \quad \textrm{as} \quad \alpha \rightarrow 0,
\end{equation*}
for all $k \in \mathbb{N}$, $\overline{\omega} \in \overline{\Omega}$ and $\omega \in \Omega$. As a consequence, we can find a sequence $(\alpha_l)_{l \in \mathbb{N}}$ converging to 0 such that
\begin{equation} \label{mollifier 2}
J_{\alpha_l} v^k(\overline{\omega}) \longrightarrow v^k(\overline{\omega}) \quad \textrm{and} \quad J_{\alpha_l} u^{\varepsilon_k}(\omega) \longrightarrow u^{\varepsilon_k}(\omega) \quad {\textrm{Leb}}_{[0,T] \times [0,\pi]} \textrm{-almost everywhere}
\end{equation}
as $l \rightarrow \infty$ for all $k \in \mathbb{N}$, $\overline{\omega} \in \overline{\Omega}$ and $\omega \in \Omega$. So~\eqref{a.e. equality in distribution v k and u eps k} follows from~\eqref{mollifier 1} and~\eqref{mollifier 2}.) 

Next, $u^{\varepsilon}$ is stochastically continuous by Theorem~4.7 in~\cite{CC2} and Lemma~B.1 in~\cite{Bally}. This and~\eqref{a.e. equality in distribution v k and u eps k} imply that $v^k$ is also stochastically continuous. By a straightforward extension of Proposition~3.21 in~\cite{Peszat} to two-parameter processes, each $v^k$ has a predictable modification $\widetilde{v}^k$. By~\eqref{pointwise convergence of v k}, we have $\widetilde{v}^{k} \longrightarrow v$ $\overline{\mathbb{P}} \otimes {\textrm{Leb}}_{[0,T] \times [0,\pi]}$-almost everywhere, so $v$ has a predictable modification as well.

Finally, the last statement is easy and we leave the details to the reader. \qed
\end{Proof}

We can now finish the proof of the weak convergence~\eqref{goal}. Indeed, the martingale problem stated in Theorem~\ref{convergence of the MPs} and satisfied by $(v, \overline{v})$ in~\eqref{Skorokhod representation 3} will allow us to identify uniquely the distribution of $(v, \overline{v})$ (from now on we may and will assume that $v$ is predictable). 

Note that the next theorem holds independently of all our previous results.

\begin{Theorem} \label{solution of MP is solution of SHE with W}
	On a filtered probability space $(\overline{\Omega}, \overline{\mathcal{F}}, \overline{\boldsymbol{F}}, \overline{\mathbb{P}})$, let $v = \{v(t,x) \mid (t,x) \in [0,T] \times [0,\pi] \}$ be an $\overline{\boldsymbol{F}}$-predictable random field and $\overline{v}$ an $\overline{\boldsymbol{F}}$-adapted càdlàg process in $H_{-r}([0,\pi])$, with $r > 1/2$. Assume that for almost all $t \in [0,T]$, $\overline{v}_t = \langle v(t, \cdot), \cdot \rangle$ as well as $\overline{v}_0=0$ holds $\overline{\mathbb{P}}$-almost surely and that	
	\begin{equation} \label{ess sup}
	\esssup_{(t,x) \in [0,T] \times [0,\pi]}  \mathbb{E}\left[ \vert v(t,x) \vert^2 \right] < \infty.
	\end{equation}
	
	In addition, assume that the pair $(v, \overline{v})$ satisfies the following martingale problem. For all $\xi \in \mathbb{R}$ and $\phi \in C_c^{\infty}((0,\pi))$, the process $(\overline{M}_t)_{t \leq T}$ defined via~\eqref{characteristics of eta} and~\eqref{def ov M} is a local $\overline{\boldsymbol{F}}$-martingale.
	
		
	Then there exists a Gaussian space--time white noise $\widetilde{W}$ on $[0,T] \times [0,\pi]$, possibly defined on a filtered extension $(\widetilde{\Omega}, \widetilde{\mathcal{F}}, \widetilde{\boldsymbol{F}}, \widetilde{\mathbb{P}})$ of $(\overline{\Omega}, \overline{\mathcal{F}}, \overline{\boldsymbol{F}}, \overline{\mathbb{P}})$ such that, with probability one, $v$ is equal in $L^2([0,T] \times [0,\pi])$ to the mild solution to the stochastic heat equation~\eqref{SHE 2} with noise $\dot{\widetilde{W}}$. Furthermore, $\overline{v}$ is indistinguishable from the modification of the latter that is continuous in $H_{-r}([0,\pi])$.
\end{Theorem}

\begin{Proof}
The proof is inspired by Lemma 2.4 in \cite{KonnoShiga}. First, Theorem II.2.42 in \cite{JJ} shows that for any $\phi \in C_c^{\infty}((0,\pi))$, the stochastic process $\langle \overline{v}, \phi \rangle$ is an $\overline{\boldsymbol{F}}$-semimartingale with first and second characteristic given by 
\begin{equation*}
t \mapsto \int_{0}^{t} \langle v(s, \cdot), \phi'' \rangle \, \textrm{d}s \quad \textrm{and} \quad 
t \mapsto \int_{0}^{t} \int_{0}^{\pi} f^2(v(s,x))  \phi^2(x) \, \textrm{d}s \, \textrm{d}x,
\end{equation*}
respectively. 
Furthermore, the third characteristic of $\langle \overline{v}, \phi \rangle $ equals 0, which implies that $\langle \overline{v}, \phi \rangle $ is continuous. As $\overline{v}_0 = 0$ $\overline{\mathbb{P}}$-almost surely, its canonical decomposition is
\begin{equation*}
\langle \overline{v}, \phi \rangle = \int_{0}^{\cdot} \langle v(s, \cdot), \phi'' \rangle \, \textrm{d}s + {\langle \overline{v}, \phi \rangle}^c,
\end{equation*}
where ${\langle \overline{v}, \phi \rangle}^c$ denotes the continuous martingale part of ${\langle \overline{v}, \phi \rangle}$. Since
\begin{equation*}
\mathbb{E} \left[\int_{0}^{T} \int_{0}^{\pi} f^2(v(s,x))  \phi^2(x) \, \textrm{d}s \, \textrm{d}x\right] \leq 
C \left( \esssup_{(s,x) \in [0,T] \times [0,\pi]} \mathbb{E}\left[ \vert v(s,x) \vert^2 \right] + 1 \right),
\end{equation*}
which is finite by assumption, the quadratic variation of ${\langle \overline{v}, \phi \rangle}^c$ is integrable, so
\begin{equation} \label{MM 1}
M_t (\phi) = \langle \overline{v}_t, \phi \rangle - \int_{0}^{t} \langle v(s, \cdot), \phi'' \rangle \, \textrm{d}s, \quad t \leq T,
\end{equation}
is a continuous square-integrable $\overline{\boldsymbol{F}}$-martingale with quadratic variation process
\begin{equation} \label{PQV 1}
t \mapsto \int_{0}^{t} \int_{0}^{\pi} f^2(v(s,x))  \phi^2(x) \, \textrm{d}s \, \textrm{d}x,
\end{equation}
for all $\phi \in C_c^{\infty}((0,\pi))$. The specifications (\ref{MM 1}) and (\ref{PQV 1}) define an orthogonal martingale measure 
$\left \{ M_t(A), \, t \in [0,T], \, A \in \mathcal{B}([0,\pi]) \right \}$
relative to $(\overline{\Omega}, \overline{\mathcal{F}}, \overline{\boldsymbol{F}}, \overline{\mathbb{P}})$, in the sense of Chapter 2 in \cite{Walsh}, with covariation measure
\begin{equation} \label{cov measure M}
Q_M(A \times B \times [s,t]) = \int_{s}^{t} \int_{A \cap B} f^2(v(r,x)) \, \textrm{d}r \, \textrm{d}x
\end{equation}
for all $A, B \in \mathcal{B}([0,\pi])$.

Now let $(\Omega', \mathcal{F}', \boldsymbol{F}', \mathbb{P}')$ be another filtered probability space on which a Gaussian space--time white noise $W'$ on $[0,T] \times [0,\pi]$ is defined. Set
\begin{equation*}
\widetilde{\Omega} = \overline{\Omega} \times \Omega', \quad 
\widetilde{\mathcal{F}} = \overline{\mathcal{F}} \otimes \mathcal{F}', \quad
\widetilde{\mathcal{F}}_t = \bigcap_{s > t}  \overline{\mathcal{F}}_s \otimes \mathcal{F}'_s, \quad
\widetilde{\mathbb{P}} = \overline{\mathbb{P}} \otimes \mathbb{P}',
\end{equation*}
and extend the random measures $M$ and $W'$ as well as the random elements $\overline{v}$ and $v$ to $\widetilde{\Omega}$ in the standard way so that on $(\widetilde{\Omega}, \widetilde{\mathcal{F}}, \widetilde{\boldsymbol{F}}, \widetilde{\mathbb{P}})$, $W'$ is independent of $(\overline{v}, v)$ and thus of $M$. In addition, on this extension, $M$ is still an orthogonal martingale measure satisfying~\eqref{MM 1} and~\eqref{cov measure M} by Lemma II.7.3 in \cite{JJ}. Define
\begin{equation*}
\begin{split}
\widetilde{W}_t(\phi) & =
\int_{0}^{t} \int_{0}^{\pi} \frac{1}{f(v(s,x))}  \mathbbm{1}_{\{ f^2(v(s,x)) \neq 0 \}} \phi(x) \, M(\textrm{d}s, \textrm{d}x) \\ & \quad \,\, +
\int_{0}^{t} \int_{0}^{\pi} \mathbbm{1}_{\{ f^2(v(s,x)) = 0 \}}  \phi(x) \, W'(\textrm{d}s, \textrm{d}x)
\end{split}
\end{equation*}
for all $t \leq T$ and $\phi \in C_c^{\infty}((0,\pi))$. As before, this defines a martingale measure 
$\{ \widetilde{W}_t(A), \, t \in [0,T], \, A \in \mathcal{B}([0,\pi]) \}$ 
relative to $(\widetilde{\Omega}, \widetilde{\mathcal{F}}, \widetilde{\boldsymbol{F}}, \widetilde{\mathbb{P}})$. 

Since $M$ and $W'$ are independent, we have from \eqref{cov measure M},
\begin{equation*}
\begin{split}
Q_{\widetilde{W}}(A \times B \times [s,t])  & =
\int_{s}^{t} \int_{A \cap B} 
\frac{1}{f^2(v(r,x))} \mathbbm{1}_{\{ f^2(v(r,x)) \neq 0 \}}  
f^2(v(r,x)) \, \textrm{d}r \, \textrm{d}x  \\ & \quad \,\, +
\int_{s}^{t} \int_{A \cap B} 
\mathbbm{1}_{\{ f^2(v(r,x)) = 0 \}}  \, \textrm{d}r \, \textrm{d}x  \\ &
= \int_{s}^{t} \int_{A \cap B}  \, \textrm{d}r \, \textrm{d}x
\end{split}
\end{equation*}
for all $A, B \in \mathcal{B}([0,\pi])$. Therefore, it follows from Proposition~2.1 in \cite{Walsh} that $\widetilde{W}$ is orthogonal and from Proposition 2.10 in \cite{Walsh} that the martingale measure $\widetilde{W}$ is a Gaussian space--time white noise on $[0,T] \times [0, \pi]$ with respect to $(\widetilde{\Omega}, \widetilde{\mathcal{F}}, \widetilde{\boldsymbol{F}}, \widetilde{\mathbb{P}})$. Moreover, we have
\begin{equation} \label{NR19}
\begin{split}
\int_{0}^{t} \int_{0}^{\pi} f(v(s,x))  \phi(x) \, \widetilde{W} (\textrm{d}s, \textrm{d}x) &  =
\int_{0}^{t} \int_{0}^{\pi} f(v(s,x))  \frac{1}{f(v(s,x))} \mathbbm{1}_{\{ f^2(v(s,x)) \neq 0 \}} \phi(x) \, M(\textrm{d}s, \textrm{d}x) \\ & \quad \,\, +
\int_{0}^{t} \int_{0}^{\pi} f(v(s,x))  
\mathbbm{1}_{\{ f^2(v(s,x)) = 0 \}}  \phi(x) \, W'(\textrm{d}s, \textrm{d}x) \\ & =
\int_{0}^{t} \int_{0}^{\pi} \mathbbm{1}_{\{ f^2(v(s,x)) \neq 0 \}}  \phi(x) \, M(\textrm{d}s, \textrm{d}x).
\end{split}
\end{equation}
Since, by~\eqref{cov measure M},
\begin{equation*}
\begin{split}
& \mathbb{E} \left[
\left(\int_{0}^{T} \int_{0}^{\pi} \mathbbm{1}_{\{ f^2(v(s,x)) \neq 0 \}}  \phi(x) \, M(\textrm{d}s, \textrm{d}x) -
\int_{0}^{T} \int_{0}^{\pi} \phi(x) \, M(\textrm{d}s, \textrm{d}x)\right)^2 \right] \\ &
\qquad = 
\mathbb{E} \left[
\left(\int_{0}^{T} \int_{0}^{\pi} \mathbbm{1}_{\{ f^2(v(s,x)) = 0 \}} \phi(x) \, M(\textrm{d}s, \textrm{d}x) \right)^2 \right] \\ &
\qquad = 
\mathbb{E} \left[ 
\int_{0}^{T} \int_{0}^{\pi} \int_{0}^{\pi} \phi(x) \mathbbm{1}_{\{ f^2(v(s,x)) = 0 \}}  
\phi(y) \mathbbm{1}_{\{ f^2(v(s,y)) = 0 \}} \, Q_M(\textrm{d}s, \textrm{d}x, \textrm{d}y) \right] \\ & \qquad  =
\mathbb{E} \left[ 
\int_{0}^{T} \int_{0}^{\pi} \phi^2(x)  \mathbbm{1}_{\{ f^2(v(s,x)) = 0 \}}  
f^2(v(s,x)) \, \textrm{d}x \, \textrm{d}s \right] = 0,
\end{split}
\end{equation*}
the $\widetilde{\boldsymbol{F}}$-martingales $t \mapsto \int_{0}^{t} \int_{0}^{\pi} \mathbbm{1}_{\{ f^2(v(s,x)) \neq 0 \}}  \phi(x) \, M(\textrm{d}s, \textrm{d}x)$ and 
$t \mapsto \int_{0}^{t} \int_{0}^{\pi} \phi(x) \, M(\textrm{d}s, \textrm{d}x)$ are indistinguishable. This implies, together with~\eqref{MM 1} and~\eqref{NR19}, that we have for any $\phi \in C_c^{\infty}((0,\pi))$,
\begin{equation} \label{weak formulation a.e.}
\begin{split}
& \int_{0}^{t} \int_{0}^{\pi} f(v(s,x)) \phi(x) \, \widetilde{W} (\textrm{d}s, \textrm{d}x) = 
M_t(\phi) = 
\langle \overline{v}_t, \phi \rangle - \int_{0}^{t} \langle v(s, \cdot), \phi'' \rangle \, \textrm{d}s, \quad t \leq T,
\end{split}
\end{equation}
$\widetilde{\mathbb{P}}$-almost surely. By assumption, the equality in~\eqref{weak formulation a.e.} holds also $\widetilde{\mathbb{P}}$-almost surely for almost all $t\leq T$ if we replace $\langle \overline{v}_t, \phi \rangle$ with $\langle v(t,\cdot), \phi \rangle$.
This and the assumption~\eqref{ess sup} imply, by the proof of Theorem 3.2 in \cite{Walsh}, that we have
\begin{equation} \label{mild formulation a.e.} 
v(t,x) = 
\int_0^t \int_{0}^{\pi} G_{t-s}(x,y) f(v(s,y)) \, \widetilde{W}(\textrm{d}s, \textrm{d}y) \quad \widetilde{\mathbb{P}}\textrm{-almost surely}
\end{equation}
for almost all $(t,x) \in [0,T] \times [0,\pi]$, i.e., $v$ satisfies the mild formulation of \eqref{SHE 2'} almost everywhere. Now let $\widetilde{v}$ be a mild solution to \eqref{SHE 2'}. Again by Theorem 3.2 in \cite{Walsh} and its proof, we can infer that $\widetilde{\mathbb{P}}$-almost surely, $v$ and $\widetilde{v}$ are equal almost everywhere and hence, in $L^2([0,T] \times [0,\pi])$. 

Finally, let $\widehat{v}$ be the continuous modification in $H_{-r}([0,\pi])$ of $\widetilde{v}$, which we obtain from Corollary~3.4 in \cite{Walsh}. By~\eqref{mild formulation a.e.}, $\widehat{v}_t = \langle v(t,\cdot), \cdot \rangle =\overline{v}_t$ $\widetilde{\mathbb{P}}$-almost surely for almost all $t\leq T$, and therefore, because $\widehat{v}$ is continuous and $\overline{v}$ càdlàg, these two processes are indistinguishable. \qed
\end{Proof}

\subsection{Necessity of the condition~\eqref{AR condition}} \label{sect:nec}

\begin{Remark} \label{Remark on f}
	Suppose that the Lipschitz function $f$ satisfies $f(0) \neq 0$. Then there must be $(t_1,x_1) \in [0,T] \times [0, \pi]$ such that 
	$\mathbb{P}( f(u(t_1,x_1)) \neq 0) > 0$, where $u$ is the mild solution to~\eqref{SHE 2}. Indeed, if we had $f(u(t,x)) = 0$ $\mathbb{P}$-almost surely for all $(t,x)$, it would imply $u = 0$ everywhere on $[0,T] \times [0, \pi]$ by equation~\eqref{SHE 2}. This in turn would imply $f(0) = 0$, which contradicts the assumption.
\end{Remark}

\begin{Theorem} \label{necessary condition}
	Assume that $f(0) \neq 0$. In the setting of Theorem~\ref{weak convergence}, if~\eqref{goal} holds, then we have~\eqref{AR condition} for all $\kappa > 0$.
\end{Theorem}

\begin{Proof}	
If~\eqref{goal} holds, we can use Skorokhod's representation theorem as in the first part of the proof of Theorem~\ref{weak convergence} and obtain for any sequence $(\varepsilon_k)_{k \in \mathbb{N}}$ converging to 0, random elements 
\begin{equation*} 
(v^{k}, \overline{v}^{k}), (v, \overline{v}) \colon (\overline{\Omega}, \overline{\mathcal{F}}, \overline{\mathbb{P}}) \longrightarrow (\Omega^*, \tau) 
\end{equation*}
on a probability space $(\overline{\Omega}, \overline{\mathcal{F}}, \overline{\mathbb{P}})$ possibly different from $(\Omega, \mathcal{F}, \mathbb{P})$ that satisfy~\eqref{Skorokhod representation 3}. Of course, we now have
\begin{equation} \label{NR26}
(v, \overline{v}) \stackrel{d}{=} (u, \overline{u}).
\end{equation}
Consider the same filtration $\overline{\boldsymbol{F}} = (\overline{\mathcal{F}}_t)_{t \leq T}$ on $\overline{\Omega}$ as in~\eqref{def filtration Om Stern}.
For fixed $\phi \in C_c^{\infty}((0,\pi))$, define the $\overline{\boldsymbol{F}}$-adapted processes
\begin{equation} \label{processes in Skorokhod rep}
\begin{split}
\overline{X}_t^k & = \langle \overline{v}^k_t, \phi \rangle - \int_{0}^{t} \int_{0}^{\pi} v^k(s,x)  \phi''(x) \, \textrm{d}s \, \textrm{d}x, \\ 
\overline{X}_t & = \langle \overline{v}_t, \phi \rangle - \int_{0}^{t} \int_{0}^{\pi} v(s,x) \phi''(x) \, \textrm{d}s \, \textrm{d}x 
\end{split}
\end{equation}
for all $k \in \mathbb{N}$ and $t \leq T$. It is straightforward to infer from~\eqref{Skorokhod representation 3} that pointwise on $\overline{\Omega}$,
\begin{equation} \label{a.s convergence in Skorokhod}
\overline{X}^k \longrightarrow \overline{X} \quad \textrm{in} \quad D([0,T], \mathbb{R}) \quad \textrm{as} \quad k \rightarrow \infty.
\end{equation}
Furthermore, by~\eqref{Skorokhod representation 3}, \eqref{NR26} and~\eqref{processes in Skorokhod rep}, $\overline{X}^k$ and $\overline{X}$ have the same distribution as the square-integrable $\boldsymbol{F}$-martingales
\begin{equation*}
t \mapsto \int_{0}^{t} \int_{0}^{\pi} 
\frac{f(u^{\varepsilon_k}(s,x))}{\sigma(\varepsilon_k)}  \phi(x) \,
L^{\varepsilon_k}(\textrm{d}s, \textrm{d}x) \quad \textrm{and} \quad 
t \mapsto \int_{0}^{t} \int_{0}^{\pi} f(u(s,x))  \phi(x) \, W (\textrm{d}s, \textrm{d}x),
\end{equation*}
respectively, and therefore, by standard arguments, we can deduce that $\overline{X}^k$ and $\overline{X}$ are $\overline{\boldsymbol{F}}$-martingales and that $\overline{X}$ is continuous.

Recall the truncation function $\varrho_h$ introduced in~\eqref{def truncation function}. Using Theorem II.2.21 in~\cite{JJ}, we can further infer that the $\overline{\boldsymbol{F}}$-semimartingale characteristics of $\overline{X}^k$ and $\overline{X}$, relative to $\varrho_h$ for a fixed but arbitrary $h > 0$, are given by 
$(\overline{B}^{k,h}, 0, \overline{\nu}^k)$ and $(0, \overline{C}, 0)$, respectively, where $\overline{\nu}^k$ is defined as in~\eqref{characteristics of eta k}, $\overline{C}$ is defined as in~\eqref{characteristics of eta} and 
\begin{equation} \label{first semimart charac k}
\overline{B}_t^{k,h} = - \int_{0}^{t} \int_{\mathbb{R}} x \mathbbm{1}_{\{ \vert x \vert > h \}} \, \overline{\nu}^k(\textrm{d} s, \textrm{d} x), \quad t \leq T.
\end{equation}
Define $\overline{X}^k(\varrho_h) = \overline{X}^k - \sum_{s \leq \cdot} \Delta \overline{X}^k_s \mathbbm{1}_{ \{ \vert \Delta \overline{X}^k_s \vert > h \}}$ for all $k \in \mathbb{N}$. Then we have, by definition of the first characteristic,
\begin{equation} \label{definition of truncation of big jumps}
\overline{X}^k(\varrho_h) = \overline{M}^{k,h} + \overline{B}^{k,h},
\end{equation}
where $\overline{M}^{k,h}$ is a local $\overline{\boldsymbol{F}}$-martingale.

Now since $\overline{X}$ is continuous, Proposition VI.2.7 in~\cite{JJ} and~\eqref{a.s convergence in Skorokhod} imply that $\overline{\omega}$-wise,
\begin{equation} \label{a.s convergence in Skorokhod 2}
\overline{X}^k(\varrho_h) \longrightarrow \overline{X} \quad \textrm{in} \quad D([0,T], \mathbb{R}) \quad \textrm{as} \quad k \rightarrow \infty.
\end{equation}
We also have  
\begin{equation}
\overline{\nu}^{k}([0,t] \times \{ \vert x \vert > a \}) \stackrel{\mathbb{\overline{P}}}{\longrightarrow} 0 \quad \textrm{as} \quad k \rightarrow \infty
\end{equation}
for any $t \leq T$ and $a > 0$ by Proposition VI.3.26 and Lemma VI.4.22 in~\cite{JJ}. Therefore, there exists a subsequence of 
$(\overline{\nu}^{k}([0,T] \times \{ \vert x \vert > h \}))_{k \in \mathbb{N}}$ converging $\overline{\mathbb{P}}$-almost surely to 0. For the sake of clarity, assume without loss of generality that this holds for the whole sequence. Applying the Cauchy--Schwarz inequality to $\overline{B}^{k,h}$ in~\eqref{first semimart charac k}, we further deduce that
\begin{equation*}
\begin{split}
\sup_{t \leq T} \big \vert \overline{B}^{k,h}_t \big \vert^2 & \leq 
\left(\int_{0}^{T} \int_{0}^{\pi} \int_{\mathbb{R}}  
\frac{f^2(v^k(t,x))}{\sigma^2(\varepsilon_k)} \phi^2(x) z^2 \, \textrm{d}t \, \textrm{d}x \, Q^{\varepsilon_k}(\textrm{d}z) \right)
\overline{\nu}^{k}([0,T] \times \{ \vert x \vert > h \}) \\ & \leq
C  \overline{\nu}^{k}([0,T] \times \{ \vert x \vert > h \})
\left(1 + \sup_{k \in \mathbb{N}} \int_{0}^{T} \int_{0}^{\pi} v^k(t,x)^2 \, \textrm{d}t \, \textrm{d}x \right)
\end{split}
\end{equation*}
and the last term converges $\overline{\mathbb{P}}$-almost surely to 0 (note that the supremum is finite because $v^k \longrightarrow v$ in $L^2([0,T] \times [0,\pi])$). This implies
\begin{equation} \label{a.s convergence in Skorokhod 3}
\overline{B}^{k,h} \longrightarrow 0 \quad \textrm{in} \quad D([0,T], \mathbb{R}) \quad \textrm{as} \quad k \rightarrow \infty
\end{equation}
$\overline{\mathbb{P}}$-almost surely. Using Proposition VI.1.23 in~\cite{JJ},~\eqref{definition of truncation of big jumps},~\eqref{a.s convergence in Skorokhod 2} and~\eqref{a.s convergence in Skorokhod 3}, we obtain
\begin{equation*}
\overline{M}^{k,h} \longrightarrow \overline{X} \quad \textrm{in} \quad D([0,T], \mathbb{R}) \quad \textrm{as} \quad k \rightarrow \infty
\end{equation*}
as well as 
\begin{equation} \label{convergence in probability}
(\overline{M}^{k,h}, -2 \overline{M}^{k,h}, ( \overline{M}^{k,h} )^2) \longrightarrow 
(\overline{X}, -2 \overline{X}, \overline{X}^2) \quad \textrm{in} \quad D([0,T], \mathbb{R}^3) \quad \textrm{as} \quad k \rightarrow \infty
\end{equation}
$\overline{\mathbb{P}}$-almost surely. Since the jumps of $\overline{M}^{k,h}$ are uniformly bounded by $h$, we can apply Proposition VI.6.13 in~\cite{JJ} on the sequence $(\overline{M}^{k,h})_{k \in \mathbb{N}}$ and then Theorem VI.6.22 (c) in~\cite{JJ} on the processes in~\eqref{convergence in probability} in order to obtain
\begin{equation*}
\left(\overline{M}^{k,h}, -2 \overline{M}^{k,h}, ( \overline{M}^{k,h} )^2, -2 \int_{0}^{\cdot} \overline{M}^{k,h}_s \, \overline{M}^{k,h} (\textrm{d}s) \right) \stackrel{\overline{\mathbb{P}}}{\longrightarrow} 
\left(\overline{X}, -2 \overline{X}, \overline{X}^2, -2 \int_{0}^{\cdot} \overline{X}_s \, \overline{X} (\textrm{d}s) \right)
\end{equation*}
in $D([0,T], \mathbb{R}^4)$ as $k \rightarrow \infty$. By definition of the quadratic variation, we can therefore deduce that
\begin{equation} \label{convergence in probability 1}
(\overline{M}^{k,h}, [\overline{M}^{k,h}, \overline{M}^{k,h}]) \stackrel{\overline{\mathbb{P}}}{\longrightarrow} (\overline{X}, \overline{C}) \quad \textrm{in} \quad D([0,T], \mathbb{R}^2) \quad \textrm{as} \quad k \rightarrow \infty.
\end{equation}
Denoting by $\overline{\mu}^k$ the jump measure of $\overline{X}^k$, we have, since $\overline{B}^{k,h}$ is continuous,
\begin{equation} \label{quadr var}
{[\overline{M}^{k,h}, \overline{M}^{k,h}]}_t = \int_{0}^{t} \int_{\mathbb{R}} x^2 \mathbbm{1}_{\{\vert x \vert \leq h\}} \, \overline{\mu}^k(\textrm{d}s , \textrm{d}x), \quad t \leq T.
\end{equation}
Now denote for any $k \in \mathbb{N}$,
\begin{equation} \label{characteristics}
\begin{split}
\widetilde{C}_t^{k} & = \int_{0}^{t} \int_{\mathbb{R}} x^2 \mathbbm{1}_{\{\vert x \vert \leq h\}} \, \overline{\nu}^k(\textrm{d}s , \textrm{d}x) \quad \textrm{and} \\
\overline{Y}_t^k & = {[\overline{M}^{k,h}, \overline{M}^{k,h}]}_t - \widetilde{C}^{k}_t = 
\int_{0}^{t} \int_{\mathbb{R}} x^2 \mathbbm{1}_{\{\vert x \vert \leq h\}} \, (\overline{\mu}^k - \overline{\nu}^k)(\textrm{d}s , \textrm{d}x).
\end{split}
\end{equation}
Then $\overline{Y}^k$ is a square-integrable $\overline{\boldsymbol{F}}$-martingale with $\vert \Delta \overline{Y}^k \vert \leq h$, and for any bounded stopping time $T$, we have, by the optional stopping theorem, $\mathbb{E}\left[ (\overline{Y}^k_T)^2 \right] \leq \mathbb{E}\left[ {[\overline{Y}^k, \overline{Y}^k]}_T \right]$. Therefore, by Lenglart's inequality (see Lemma I.3.30 in~\cite{JJ}), we obtain for all $\delta >0$ and $\eta > 0$,
\begin{equation} \label{Lenglart}
\begin{split}
\mathbb{P}\left( \sup_{s \leq t} \vert \overline{Y}^k_s \vert^2 \geq \delta \right) & \leq \frac{1}{\delta} 
\left( \eta + \mathbb{E}\left[ \sup_{s \leq t} \Delta {[\overline{Y}^k, \overline{Y}^k]}_s \right] \right) + \mathbb{P}({[\overline{Y}^k, \overline{Y}^k]}_t \geq \eta) \\ & \leq
2 \frac{\eta}{\delta} + \left(\frac{h}{\delta} + 1\right)\mathbb{P}({[\overline{Y}^k, \overline{Y}^k]}_t \geq \eta).
\end{split}
\end{equation}
By~\eqref{quadr var}, we have
\begin{equation*}
{[\overline{Y}^k, \overline{Y}^k]}_t = \int_0^t \int_{\mathbb{R}} x^4 \mathbbm{1}_{\{\vert x \vert \leq h\}} \, \overline{\mu}^k(\textrm{d}s , \textrm{d}x) \leq
\left(\sup_{s \leq t} \big \vert \Delta {[\overline{M}^{k,h}, \overline{M}^{k,h}]}_s  \big \vert \right) {[\overline{M}^{k,h}, \overline{M}^{k,h}]}_t.
\end{equation*}
Moreover, because ${[\overline{M}^{k,h}, \overline{M}^{k,h}]}_t \stackrel{\overline{\mathbb{P}}}{\longrightarrow} \overline{C}_t$ by~\eqref{convergence in probability 1} and $\sup_{s \leq t} \vert \Delta {[\overline{M}^{k,h}, \overline{M}^{k,h}]}_s \vert \stackrel{\overline{\mathbb{P}}}{\longrightarrow} 0$ by Proposition VI.3.26 (iii) in~\cite{JJ}, we deduce from the inequality above that ${[\overline{Y}^k, \overline{Y}^k]}_t \stackrel{\overline{\mathbb{P}}}{\longrightarrow} 0$ and, by~\eqref{Lenglart}, that
\begin{equation} \label{convergence in probability 2}
\sup_{s \leq t} \vert \overline{Y}^k_s \vert \stackrel{\overline{\mathbb{P}}}{\longrightarrow} 0 \quad \textrm{as} \quad k \rightarrow \infty
\end{equation}
for all $t \leq T$. Finally, combine~\eqref{convergence in probability 1},~\eqref{characteristics} and~\eqref{convergence in probability 2} to see that
\begin{equation} \label{convergence in probability 3}
\widetilde{C}_t^k = \int_{0}^{t} \int_{\mathbb{R}} x^2 \mathbbm{1}_{\{\vert x \vert \leq h\}} \, \overline{\nu}^k(\textrm{d}s , \textrm{d}x) \stackrel{\overline{\mathbb{P}}}{\longrightarrow} \overline{C}_t \quad \textrm{as} \quad k \rightarrow \infty
\end{equation}
for all $t \leq T$ and $h > 0$. Taking a subsequence if necessary, we will from now on assume that the convergence in~\eqref{convergence in probability 3} holds even $\overline{\mathbb{P}}$-almost surely.

Recall now the definition of $\Sigma_h^k(s,x)$ in~\eqref{def I and Sigma} and that, because $v^k \longrightarrow v$ in $L^2([0,T] \times [0,\pi])$, we have $\int_{0}^{t} \int_{0}^{\pi}
\vert f^2(v^k(s,x)) - f^2(v(s,x)) \vert  \phi^2(x) \,\textrm{d}s \, \textrm{d}x \longrightarrow 0$ as $k \rightarrow \infty$; see the calculations in~\eqref{NR10} and~\eqref{NR11}. Together with~\eqref{convergence in probability 3} this implies $\overline{\mathbb{P}}$-almost surely, 
\begin{equation} \label{convergence}
\int_{0}^{T} \int_{0}^{\pi}
f^2(v(s,x)) \phi^2(x) \Sigma_h^k(s,x) \,\textrm{d}s \, \textrm{d}x \longrightarrow 0 \quad \textrm{as} \quad 
k \rightarrow \infty
\end{equation}
for all $h>0$ and $\phi \in C_c^{\infty}((0,\pi))$ by a similar calculation as in~\eqref{NR1} (note that the first inequality there becomes an equality if $\vert \cdot \vert$ is replaced by $(\cdot)$ throughout).

Now on the set $\{ \vert f(v^k(s,x)) \phi(x) \vert \geq \delta \}$, where $\delta >0$, we have 
\begin{equation*}
\mathbbm{1}_{ \left \{ \vert z \vert \geq \left(h/\vert f(v^k(s,x)) \phi(x) \vert\right) \sigma(\varepsilon_k) \right \}} \geq 
\mathbbm{1}_{ \left \{ \vert z \vert \geq \left(h/\delta\right) \sigma(\varepsilon_k) \right \}},
\end{equation*}
and thus
\begin{equation*}
\Sigma_h^k(s,x) \geq \frac{1}{\sigma^2(\varepsilon_k)} 
\int_{\mathbb{R}} z^2  
\mathbbm{1}_{ \left \{ \vert z \vert \geq \left(h/\delta\right)  \sigma(\varepsilon_k) \right \}} \, Q^{\varepsilon_k}(\textrm{d}z).
\end{equation*}
Therefore, as a consequence of~\eqref{convergence}, we obtain $\overline{\mathbb{P}}$-almost surely,
\begin{equation} \label{NR17}
\begin{split}
& \frac{1}{\sigma^2(\varepsilon_k)}  
\int_{\mathbb{R}} z^2  
\mathbbm{1}_{ \left \{ \vert z \vert \geq \left(h/\delta\right) \sigma(\varepsilon_k) \right \}} \, Q^{\varepsilon_k}(\textrm{d}z) \\ & \qquad \times
\int_{0}^{T} \int_{0}^{\pi}
f^2(v(s,x))  \phi^2(x)  \mathbbm{1}_{\{ \vert f(v^k(s,x)) \phi(x) \vert \geq \delta, \, \vert f(v(s,x)) \phi(x) \vert > \delta \}} \,\textrm{d}s \, \textrm{d}x \longrightarrow 0
\end{split}
\end{equation}
as $k \rightarrow \infty$ for all $h>0$, $\delta>0$ and $\phi \in C_c^{\infty}((0,\pi))$.

We have seen in~\eqref{pointwise convergence of v k} that we can assume (perhaps for a subsequence) that 
\begin{equation*}
v^{k} \longrightarrow v \quad \textrm{as} \quad k \rightarrow \infty \quad \overline{\mathbb{P}} \otimes {\textrm{Leb}}_{[0,T] \times [0,\pi]}\textrm{-almost everywhere},
\end{equation*}
which implies, by dominated convergence and continuity of $f$,
\begin{equation} \label{NR20}
\begin{split}
& \mathbb{E}\left[
\int_{0}^{T} \int_{0}^{\pi}
f^2(v(s,x)) \phi^2(x) 
\mathbbm{1}_{\{ \vert f(v^k(s,x))  \phi(x) \vert \geq \delta, \, \vert f(v(s,x)) \phi(x) \vert > \delta \}} \,\textrm{d}s \, \textrm{d}x \right] \\ & \qquad \longrightarrow
\mathbb{E}\left[
\int_{0}^{T} \int_{0}^{\pi}
f^2(v(s,x)) \phi^2(x) 
\mathbbm{1}_{\{ \vert f(v(s,x)) \phi(x) \vert > \delta \}} \,\textrm{d}s \, \textrm{d}x \right] \quad \textrm{as} \quad k \rightarrow \infty.
\end{split}
\end{equation}
So from~\eqref{Skorokhod representation 3},~\eqref{NR17} and~\eqref{NR20}, we deduce that
\begin{equation} \label{NR21}
\begin{split}
& \frac{1}{\sigma^2(\varepsilon_k)}  
\int_{\mathbb{R}} z^2  
\mathbbm{1}_{ \left \{ \vert z \vert \geq \left(h/\delta\right) \sigma(\varepsilon_k) \right \}} \, Q^{\varepsilon_k}(\textrm{d}z) \\ & \qquad \times
\mathbb{E} \left[\int_{0}^{T} \int_{0}^{\pi}
f^2(u(s,x))  \phi^2(x)  \mathbbm{1}_{\{ \vert f(u(s,x))  \phi(x) \vert > \delta \}} \,\textrm{d}s \, \textrm{d}x\right] \longrightarrow 0 \quad \textrm{as} \quad 
k \rightarrow \infty
\end{split}
\end{equation}
for all $h>0$, $\delta>0$ and $\phi \in C_c^{\infty}((0,\pi))$. Moreover,
\begin{equation} \label{NR22}
\begin{split}
& \left( \frac{1}{\sigma^2(\varepsilon_k)} 
\int_{\mathbb{R}} z^2  
\mathbbm{1}_{ \left \{ \vert z \vert \geq \left(h/\delta\right) \sigma(\varepsilon_k) \right \}} \, Q^{\varepsilon_k}(\textrm{d}z) \right)
\mathbb{E} \left[\int_{0}^{T} \int_{0}^{\pi}
f^2(u(s,x))  \phi^2(x) \, \textrm{d}s \, \textrm{d}x\right] 
\\ & \qquad \leq \left( \frac{1}{\sigma^2(\varepsilon_k)} 
\int_{\mathbb{R}} z^2  
\mathbbm{1}_{ \left \{ \vert z \vert \geq \left(h/\delta\right)  \sigma(\varepsilon_k) \right \}} \, Q^{\varepsilon_k}(\textrm{d}z) \right) \\ & \qquad \quad \,\, \times
\mathbb{E} \left[\int_{0}^{T} \int_{0}^{\pi}
f^2(u(s,x))  \phi^2(x)  \mathbbm{1}_{\{ \vert f(u(s,x))  \phi(x) \vert > \delta \}} \,\textrm{d}s \, \textrm{d}x\right] + T \pi \delta^2.
\end{split}
\end{equation}
So if we choose $h = \kappa \delta$ with $\kappa > 0 $ arbitrary, then by~\eqref{NR21}, the first term on the right-hand side of~\eqref{NR22} converges to 0 as $k \rightarrow \infty$ for all $\kappa >0$ and $\delta > 0$. The second term does not depend on $k$ nor $h$ and converges to 0 as $\delta \rightarrow 0$. This implies
\begin{equation} \label{NR5}
\left( \frac{1}{\sigma^2(\varepsilon_k)} 
\int_{\mathbb{R}} z^2  
\mathbbm{1}_{ \left \{ \vert z \vert \geq \kappa \sigma(\varepsilon_k) \right \}} \, Q^{\varepsilon_k}(\textrm{d}z) \right)
\mathbb{E} \left[\int_{0}^{T} \int_{0}^{\pi}
f^2(u(s,x))  \phi^2(x) \, \textrm{d}s \, \textrm{d}x\right] \longrightarrow 0 
\end{equation}
as $k \rightarrow \infty$ for all $\kappa > 0 $. Since $f(0) \neq 0$, there exists $(t_1,x_1) \in [0,T] \times [0, \pi]$ such that 
$\mathbb{E} [ f^2(u(t_1,x_1)) ]  > 0$ by Remark~\ref{Remark on f}. Moreover, the mild solution $u$ is continuous in $L^2(\Omega, \mathcal{F},\mathbb{P})$, which follows from the proof of Corollary~3.4 in~\cite{Walsh}. We can thus infer that the expectation in~\eqref{NR5} is not 0 and we obtain
\begin{equation*}
\frac{1}{\sigma^2(\varepsilon_k)}  
\int_{\mathbb{R}} z^2  
\mathbbm{1}_{ \left \{ \vert z \vert \geq \kappa  \sigma(\varepsilon_k) \right \}} \, Q^{\varepsilon_k}(\textrm{d}z) \longrightarrow 0 \quad \textrm{as} \quad k \rightarrow \infty
\end{equation*}
for all $\kappa > 0$, which is exactly~\eqref{AR condition}. \qed

\end{Proof}

\subsection*{Acknowledgements}
TD cordially thanks Claudia Klüppelberg for inspiring discussions and valuable advice as well as the Chair of Probabilities of the École Polytechnique Fédérale de Lausanne for its hospitality during his visit.
TD's research is partially supported by the Deutsche Forschungsgemeinschaft, project number KL 1041/7-1.

\addcontentsline{toc}{section}{References}
\bibliographystyle{plainnat}
\bibliography{bib-Normal_Approximation_revised}

\end{document}